**AUTONOMOUS STRUCTURAL DISCOVERY OF BIOPROCESS DYNAMICS AND KINETICS VIA LLM-DRIVEN AGENTS**

Erdal Aydin[1,*], Yaman Arkun[2**]

[1]School of Chemical Engineering, Aalto University, Espoo 02150, Finland

[2]Department of Chemical and Biological Engineering, Koc University, 34450, Istanbul, Türkiye

[*]erdal.aydin@aalto.fi

[**]yarkun@ku.edu.tr

**Abstract**

High-fidelity digital twins in continuous biomanufacturing are often bottlenecked by incomplete kinetic models, while standard data-driven algorithms (e.g., SINDy, PINNs) lack interpretability or impose rigid physical assumptions. This paper introduces a hybrid closed-loop architecture for autonomous kinetic state-space discovery driven by a Large Language Model (LLM) agent. The framework employs a sequential "isolation-peeling" methodology across orthogonal operational regimes. This algorithmic peeling mathematically decouples entangled multiplicative (inhibitory) and additive (catalytic) biological rate law. Grounded in a Two-Level Retrieval-Augmented Generation (RAG) library, the framework utilizes Moving Horizon Estimation (MHE) and logarithmic decomposition to decouple complex cell growth kinetics into orthogonal multiplicative (inhibitory) and additive (regulatory) residual channels. When statistical degeneracy occurs an Active Design of Experiments (DoE) controller dynamically triggers targeted perturbations to decorrelate competing model hypotheses. Evaluated on a continuous stirred-tank reactor exhibiting unmodeled Haldane inhibition, Luong product toxicity, and secondary metabolite autocatalysis, the framework autonomously expanded the state space and recovered ground-truth parameters with under 6% error without a priori kinetic knowledge.

## 1. INTRODUCTION

Continuous biomanufacturing has emerged as a strategic priority for the biopharmaceutical industry, promising superior product consistency, a smaller physical footprint, and faster responsiveness to fluctuating drug demand relative to legacy batch processing (Schmidt et al., 2022). Regulatory agencies have actively encouraged this transition through Process Analytical Technology (PAT) initiatives that mandate real-time, in-line characterization of critical process parameters. Realizing the full economic and quality benefits of continuous operation, however, is contingent on a digital twin capable of tracking and predicting bioreactor behavior with high fidelity across its entire operating envelope, a task that becomes acutely difficult whenever the twin's underlying kinetic model is structurally incomplete.

The deployment of high-fidelity digital twins and advanced model predictive control strategies is central to this transition, and in turn depends on the accurate representation of process kinetics. Unlike classical chemical reactors, biological Continuous Stirred-Tank Reactors (CSTRs) exhibit highly non-linear, time-varying behaviors, including substrate inhibition, secondary metabolite toxicity, and complex competitive interactions. Traditional first-principles modeling requires exhaustive *a priori* knowledge, while purely data-driven approaches often fail to generalize outside their training domain and inherently lack mechanistic interpretability.

Consequently, hybrid modeling and system identification have emerged as critical research frontiers. State estimation techniques, such as the Extended Kalman Filter (EKF) and Moving Horizon Estimation (MHE), excel at tracking unmeasured states and estimating time-varying parameters (Rawlings & Mayne, 2009). However, these methods operate under the strict assumption that the underlying mathematical structure of the model is correct. When a system encounters emergent, unmodeled phenomena—such as unexpected structural toxicity, secondary metabolite accumulation, or shifting inhibition modalities—passive parameter estimation becomes entirely insufficient. The fundamental challenge thus shifts from parameter tuning to dynamic model state-space expansion and structural kinetic discovery.

Recent advancements have explored symbolic regression and sparse identification algorithms (e.g., SINDy) to uncover governing equations directly from data (Brunton et al., 2016). While mathematically powerful, these methods often prove inadequate for rigorous bioprocess modeling. Because SINDy promotes sparsity across a library of candidate monomials, it frequently yields high-dimensional polynomial approximations that lack direct biological meaning. In a bioprocess context, practitioners require the identification of specific, physically bounded constants—such as half-saturation coefficients ($K_s$) or inhibition thresholds ($K_i$)—rather than abstract mathematical weights (Chang et al., 2008). Divorced from structural mechanistic constraints, standard symbolic regression often produces models that fundamentally drift from biological reality. Building on this need for physically interpretable structures, del Rio-Chanona and coworkers have proposed automated discovery frameworks that pair genetic programming with information-theoretic model selection to recover catalytic and bioprocess kinetic rate laws (de Carvalho Servia et al., 2024; Mowbray et al., 2022); however, these frameworks still require the candidate mechanism space to be defined in advance and do not autonomously expand the state-space dimensionality when a governing species is entirely absent from the model.

Similarly, the application of deep learning—specifically Neural Networks (NNs) and Physics-Informed Neural Networks (PINNs)—introduces significant barriers to mechanistic transparency. While NNs can interpolate highly non-linear datasets with extreme precision, they function as opaque "black-box" structures whose internal biases offer no actionable physical insight (LeCun et al., 2015). PINNs attempt to bridge this gap by embedding differential equations directly into the loss function (Raissi et al., 2019); yet, they remain fundamentally constrained by structural rigidity. PINNs are designed to solve systems where the governing physics are already well-defined. If the underlying bioprocess model lacks a crucial mass balance or misidentifies a kinetic rate law, PINNs will effectively force the empirical data to adhere to an incorrect structural hypothesis.

In parallel with these developments in numerical system identification, Large Language Models (LLMs) have rapidly matured from static text generators into autonomous, tool-using scientific agents capable of formulating hypotheses, executing external computations, and revising their beliefs based on empirical feedback. Systems such as Coscientist have shown that an LLM agent, given access to laboratory automation application programming interfaces (APIs), can independently plan and execute multi-step chemical experiments (Boiko et al., 2023), while ChemCrow augments an LLM with domain-specific computational tools to solve organic synthesis and drug discovery tasks that exceed the model's native reasoning ability (M. Bran et al., 2024). Central to this emerging agentic paradigm is Retrieval-Augmented Generation (RAG), which grounds an LLM's reasoning in an external, structured corpus of domain knowledge rather than relying solely on parametric memory (Lewis et al., 2020). Unlike the black-box function approximators discussed above, an LLM agent equipped with a curated library of mechanistic kinetic forms can reason explicitly over candidate biological hypotheses, select among them using interpretable statistical criteria, and, critically, judge when the existing state-space representation itself is inadequate and must be expanded. This capacity for structural, rather than purely parametric, reasoning is precisely what is missing from EKF/MHE, SINDy, and PINN-based approaches, and motivates our use of an LLM-driven agent as the core diagnostic engine of the discovery framework proposed here. This agentic paradigm is now reaching process systems engineering directly: multi-agent LLM frameworks have recently been shown to automate operational decision-making, flowsheet synthesis, and simulation workflows in chemical process contexts (Pajak et al., 2025), lending further precedent to the use of an LLM agent as an autonomous experimental and diagnostic collaborator in a regulated bioprocess setting.

Therefore, a critical bottleneck remains in the inability of current diagnostic architectures to autonomously correct structural model deficits. To fully realize autonomous continuous biomanufacturing, a framework is required that can both insert missing physical states (ODE generation) and subsequently resolve the true, non-linear kinetic topologies governing those states.

The contributions of this work are as follows: First, we introduce a hybrid discovery architecture in which an LLM agent, grounded in a domain-specific Two-Level Retrieval-Augmented Generation (RAG) library, autonomously selects, parameterizes, and validates candidate kinetic structures against noisy bioreactor telemetry. Second, we develop a Moving Horizon Estimation and logarithmic-decomposition protocol that isolates the specific growth rate into orthogonal multiplicative (inhibitory) and additive (regulatory) residual channels, enabling sequential, interference-free structural identification. Third, we design an orthogonal Design-of-Experiments protocol, coupled with an Active DoE controller, that allows the agent to autonomously trigger targeted substrate perturbations whenever competing kinetic hypotheses become statistically degenerate ($\Delta AIC < 2.0$). Finally, we demonstrate, on a synthetic CSTR exhibiting unmodeled Haldane inhibition, Luong toxicity, and autocatalytic secondary-metabolite kinetics, that the resulting framework recovers the full five-dimensional state-space and the ground-truth kinetic parameters to within 6% error, without any *a priori* knowledge of the missing physics.

The remainder of this paper is organized as follows. Section 2 describes the synthetic bioprocess plant, its multi-modal kinetic manifold, and the steady-state parameter estimation procedure used to anchor known stoichiometric and Monod constants. Section 3 details the autonomous discovery architecture, including the Moving Horizon Estimation protocol, the orthogonal Design-of-Experiments regimes, the Two-Level RAG library, and the structural identification and Active DoE intervention logic. Section 4 presents the results of the discovery tournaments and the global parameter refinement. Section 5 discusses the limitations of the current hybrid offline/online architecture, and Section 6 outlines directions for extending the framework toward open-world, fully autonomous discovery.

## 2. DESCRIPTION OF THE SYSTEM

### 2. 1. Plant Description and Modeling Assumptions

To evaluate the discovery architecture, we consider a non-linear CSTR where the biological growth process is governed by a complex, multi-modal kinetic manifold with the following characteristics

- **Substrate Inhibition:** Modeled via Haldane kinetics ($\mu_H$) to account for metabolic overloading at high substrate concentrations.
- **Product Toxicity:** Represented by a Luong-type product inhibition function ($\mu_L$) that captures the deleterious and growth-limiting impact of the accumulating primary product ($P$).
- **Autocatalysis:** A low-affinity term ($\mu_{auto}$) accounting for auxiliary metabolic pathways triggered by the accumulation of a specific secondary precursor ($Ps$) that drives the autocatalytic growth loop.

The specific growth rate of the true plant is defined by a multiplicative manifold representing inhibition, augmented by additive metabolic pathways:

$$\mu = \left( \underbrace{\frac{\mu_{max} S}{K_s + S}}_{\text{Monod Growth } (\mu_{base})} \cdot \underbrace{\frac{1}{1 + \frac{S^2}{K_I (15 + S)}}}_{\text{Haldane Substrate Inhibition} (\mu_H)} \right) \cdot \underbrace{\left(1 - \frac{P}{P_m}\right)}_{\mu_L} + \underbrace{\frac{k_{auto} Ps}{K_p + Ps}}_{\mu_{auto}}$$

or

$$\mu = \underbrace{\left(\mu_{base} \cdot \mu_H(S) \cdot \mu_L(P)\right)}_{\text{Multiplicative Inhibition Path}} + \underbrace{\mu_{auto}}_{\text{Additive Rescue Path}} \tag{1}$$

The term $\mu_{base}$ represents the fundamental Monod Growth term that describes the saturation of enzymatic activity at low substrate concentrations. While the Monod term wants growth to increase, substrate inhibition of Haldane type $\mu_H$ forces it to slow down at high substrate levels. Physically, when substrate concentration is too high, it binds to the enzymes in a way that blocks the active site or disrupts the internal osmotic balance of the cell. Product inhibition ($\mu_P$) acts as a process feedback brake. As the cells produce $P$ (the product), the growth decreases and a rapid collapse in metabolic efficiency occurs once a critical toxicity threshold $P_m$, is reached.

The secondary metabolite ($P_s$) acts as a metabolic "turbo" for the system. As the metabolite accumulates via the cell's own biosynthetic activity, it triggers secondary pathways that provide an additive boost $\mu_{auto}$ to the total specific growth rate ($\mu$). This creates a highly nonlinear, autocatalytic positive feedback loop: cellular growth drives the accumulation of $P_s$, and the elevated $P_s$ concentration, in turn, further accelerates growth. The "True plant" dynamics is governed by the following system of non-linear differential equations which are simulated to get the "real" data:

$$\frac{dX}{dt} = (\mu - D)X \tag{2}$$

$$\frac{dS}{dt} = D(S_{in} - S) - \frac{\mu}{Y_{x/s}} X \tag{3}$$

$$\frac{dP}{dt} = Y_{p/x} \mu X - D(P - P_{in}) \tag{4}$$

$$\frac{dPs}{dt} = Y_{Ps/X}\mu X - D \cdot Ps + Ps_{in} \tag{5}$$

The dynamic profile of the secondary metabolite ($P_s$) reflects a cellular shift from primary biomass accumulation to secondary metabolism. The specific production yield ($Y_{P_s/X}$) is governed by a "Gated Hill" architecture, which mathematically defines this transition through the product of two distinct mechanistic terms:

$$Y_{P_s/X} = \left(k_{base} + k_{max}\frac{1}{1 + (S/K)^n}\right)\left(\frac{S}{K_{gate} + S}\right) \tag{6}$$

The first bracketed term denotes the transcriptional starvation switch (the inhibitory Hill function). Secondary metabolic pathways are highly resource-intensive and are typically repressed during periods of substrate abundance. Under normal operational baselines, the $(S/K)^n$ term dominates the denominator, effectively collapsing the yield to its minimal basal rate ($k_{base}$).

The second term, $S/(K_{gate} + S)$, acts as a strict physical boundary condition. While the biological switch dictates that cells "attempt" to maximize secondary production under severe starvation, biosynthesis intrinsically requires a carbon and energy source. As the substrate concentration approaches absolute zero, this thermodynamic gate overrides the regulatory switch, forcing the overall production yield back to zero. The constant $K_{gate}$ operates entirely at the microscopic physical limit and is structurally decoupled from the macroscopic threshold $K$.

A plot of the specific secondary metabolite yield ($Y_{P_s/X}$) against substrate concentration ($S$) is shown in Figure 1. The region to the right of the starvation threshold ($S > K$) represents the nutrient-rich state where the cells prioritize consuming the substrate solely to build biomass ($X$) and primary product ($P$). The region to the left of the threshold ($S < K$) represents the nutrient-limited state where primary growth halts, and the cells begin producing their peculiar secondary metabolites ($P_s$).

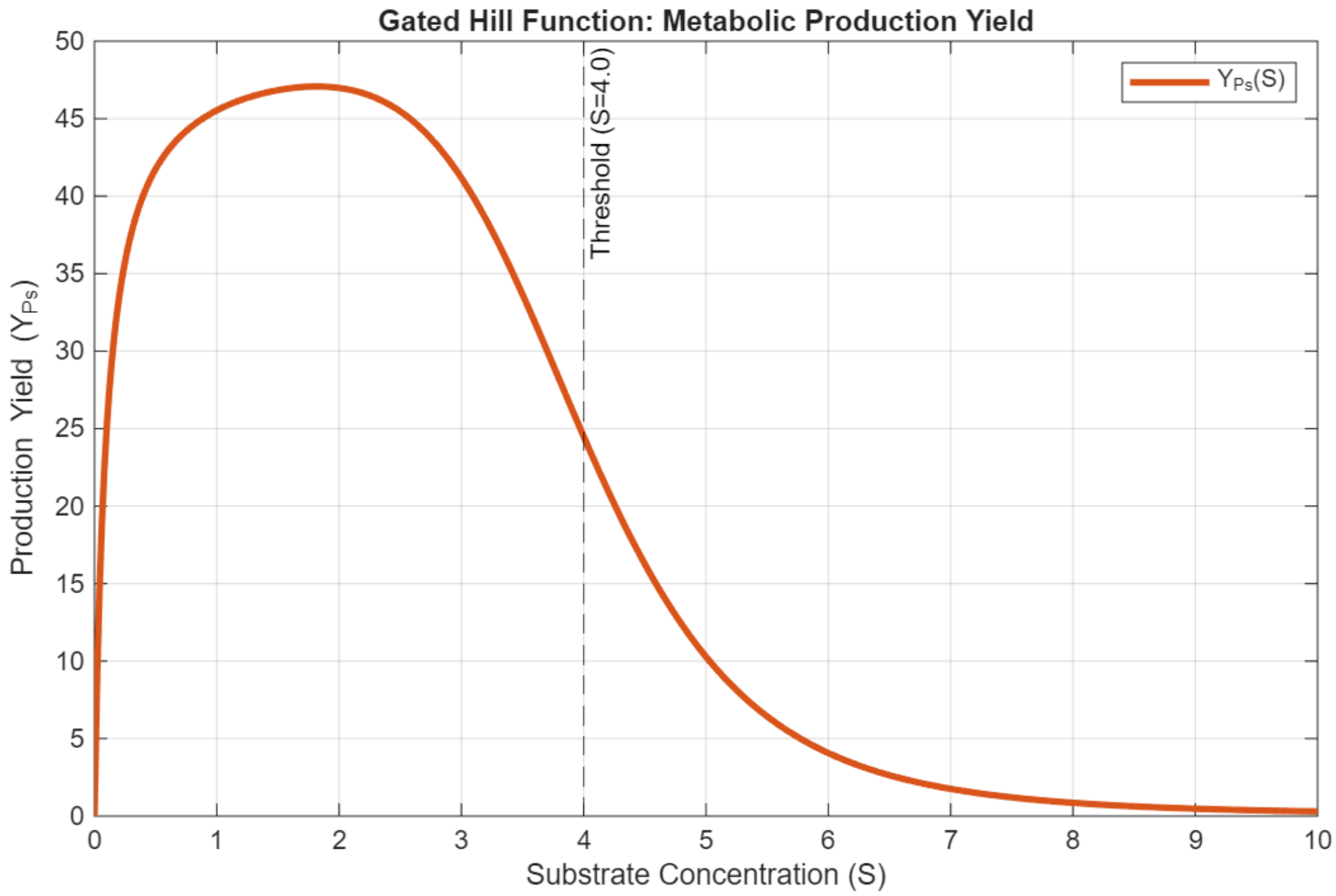

**Figure 1.** Product Yield $Y_{P_s/X}$ for the secondary metabolite $P_s$.

In this system, $\mu$ is treated as an unobserved latent variable. The states are concentrations of the biomass X, the substrate, S, the main product, P, and the secondary precursor Ps. Inputs are $D, S_{in}, P_{in}$ $and$ $Ps_{in}$. Depending on the experimental protocol, they are used as manipulated variables by the feedback controller to regulate the states or as pulse inputs.

The total specific growth rate ($\mu_{total}$) that governs the true dynamic trajectory of the biomass is the direct summation of the invariant fundamental carbon pathway, the mass-dependent toxicity suppression ($\mu_{multip}$), and the feed-scaled secondary additive boost ($\mu_{additive}$). By applying the steady-state mass balance constraints—where biomass and product concentrations resolve to algebraic dependencies of $S$—the time-invariant equilibrium manifold is isolated and plotted as shown in Figure 2. The curve is generated under the strict assumption that the feed stream is sterile and product-free ( $X_{in} = 0$, $P_{in} = 0$, and $P_{s,in} = 0$) and the carbon load $S_{in}$ is free to change.

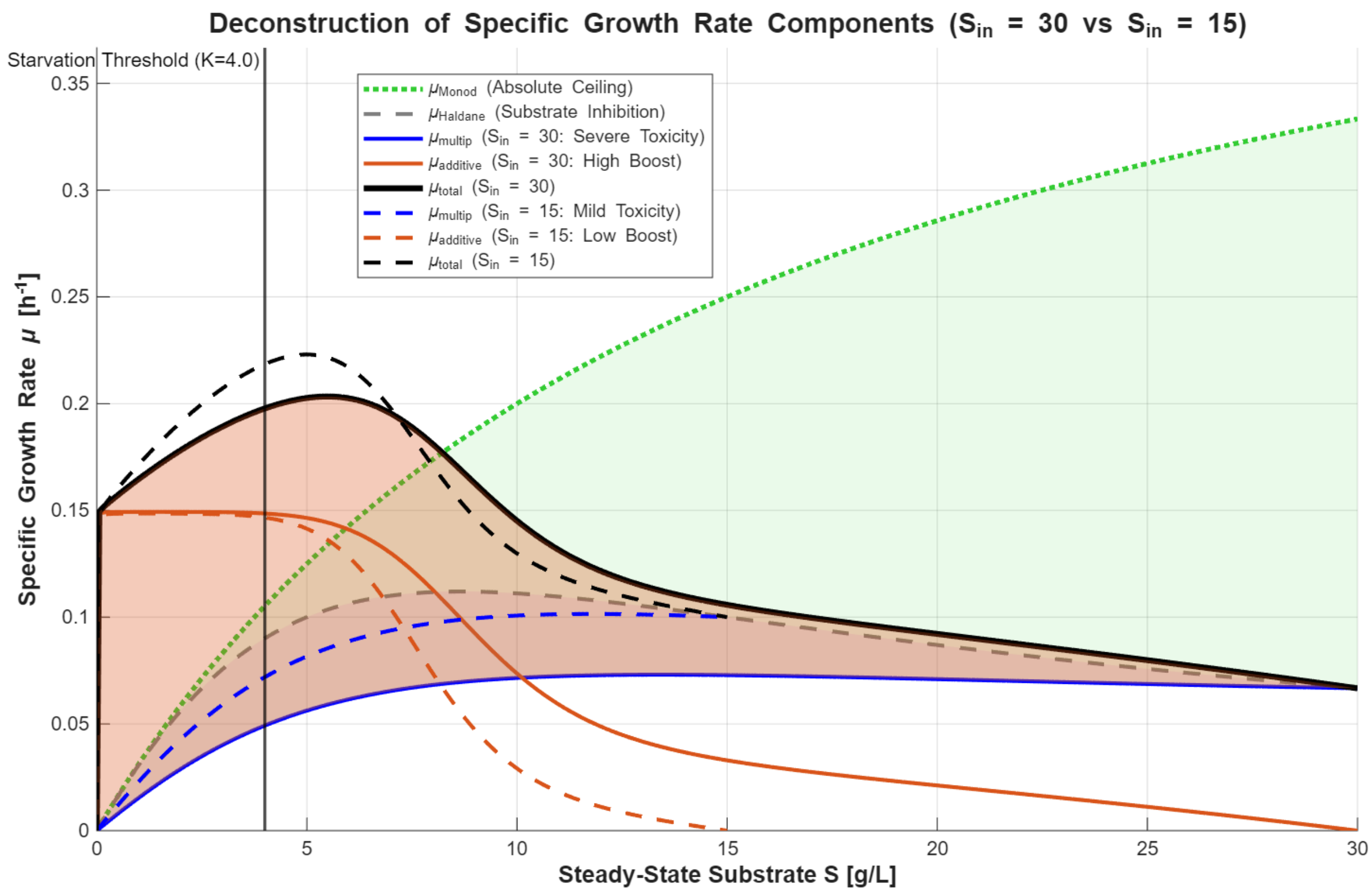

**Figure 2.** Steady state specific growth rate with its individual constituents**.**

During the dynamic model discovery phase, the bioreactor is not operated under a static feed regime. To generate sufficient transient variance for the discovery algorithm, a closed-loop feedback controller manipulates the inlet substrate concentration ($S_{in}$) to track programmed shifts in the substrate setpoint ($S_{sp}$). By superimposing two distinct operational boundaries ($S_{in}$ = 15.0 g/L and $S_{in}$ = 1.0 g/L), the effect of $S_{in}$ on the steady-state kinetic manifolds are visualized in Figure 2. As illustrated, while the differing mass fluxes shift the equilibrium positions of the curves, the fundamental topological shape of the kinetic manifolds remains strictly invariant, confirming that the underlying biological rate laws are conserved regardless of the external feed load.

### 2.2.1. The Invariant Primary Kinetics

The theoretical ceiling of cellular proliferation, defined by classic uninhibited Monod kinetics ($\mu_{Monod}$), and the primary substrate inhibition, modeled via Haldane kinetics ($\mu_{Haldane}$), represent fundamental microscopic interactions. Because these mechanisms are governed solely by the immediate, localized substrate concentration ($S$) encountered by the cells, their mathematical topologies remain perfectly rigid and completely invariant to the total carbon load pumped into the reactor.

### 2.2.2 The Mass-Dependent Toxicity Penalty

In contrast to the primary kinetics, the sequential penalties applied to the growth baseline are strictly endogenous and scale directly with $S_{in}$. Under the high-carbon feed ($S_{in} = 50.0$ g/L), macroscopic mass conservation dictates that the reactor converts the massive influx of carbon into a proportionally massive concentration of primary product ($P$). This accumulation triggers severe Luong toxicity, creating a highly suppressed primary growth regime ($\mu_{multip}$). Conversely, under the lower carbon load ($S_{in} = 25.0$ g/L), absolute product accumulation is restricted, resulting in a substantially delayed and milder toxicity penalty.

### 2.2.3. The Feed-Scaled Autocatalytic Boost

Operating in parallel to this bottlenecked primary pathway is the autocatalytic positive feedback loop triggered during substrate starvation. The independent metabolic boost ($\mu_{additive}$) provided by the accumulation of the secondary metabolite ($P_s$) is similarly governed by macroscopic scaling. The boost is highly pronounced at $S_{in} = 50.0$ g/L due to the larger absolute mass of available carbon being diverted into the secondary metabolic pathways once the starvation threshold is crossed. At $S_{in} = 25.0$ g/L, the diminished carbon availability results in a proportionally smaller secondary boost.

The discovery agent is tasked with reconstructing the functional form of $\mu$ through the observation of noisy state measurements, with the constraint that the agent has no *a-priori* knowledge of the inhibitory mechanisms or the autocatalytic parameters. The "True Plant" operates on a 4-dimensional state vector with dynamics defined by equations (2)-(5). However, by design, the state space of base model is structurally constrained by three states (S, X and P) excluding the secondary metabolite $P_s$. It also relies on an elementary Monod growth assumption that is insufficient to capture the complex, multi-variable kinetic landscape of the true plant. Thus, the discovery engine is faced with a plant model mismatch due to missing physics (inhibitory kinetics) and missing structure due to the dynamics of the secondary metabolite $P_s$. Ultimately, we do not propose this AI framework as a replacement for first-principles modeling. Rather, it serves as an autonomous structural diagnostic tool that detects where initial first-principles models are incomplete, allowing us to iteratively expand the state space of our kinetic models until they fully match the mechanistic complexity of the true plant.

### 2.3 Parameter Estimation via Steady-State Baselines

To ensure absolute structural identifiability during the dynamic model discovery phase, the foundational stoichiometric and kinetic parameters must be mathematically decoupled from the unknown transient dynamics. Attempting to identify all system parameters simultaneously within a nonlinear dynamic framework frequently results in an ill-posed optimization problem, leading to parameter collinearity and the emergence of structural decoys. To circumvent this, fundamental mass balances, growth affinities, and metabolic yield boundaries can be independently estimated using offline and steady-state experimental data, and subsequently fixed as rigid mathematical anchors.

#### 2.3.1. Offline Estimation of Primary Stoichiometric Yields ($Y_{X/S}$ and $Y_{P/X}$)

The biomass yield on substrate ($Y_{X/S}$) and the primary product yield on biomass ($Y_{P/X}$) represent intrinsic metabolic conversion efficiencies that remain largely constant during balanced primary growth. These parameters are calculated via linear regression of offline batch or steady-state analytical measurements. For a continuous system at steady state, the relationships decouple algebraically from the time domain:

$$X = Y_{X/S}(S_{in} - S) \tag{7}$$

$$P = Y_{P/X}X \tag{8}$$

By plotting the steady-state biomass concentration ($X$) against the total substrate consumed ($S_{in} - S$), and the primary product ($P$) against the biomass, the specific yield constants can be empirically determined and fixed during the structural discovery phase. Locking these fundamental yield constants, establishes a rigid, mass-conserving backbone for the digital twin. This guarantees that the dynamic model cannot falsely adjust basic stoichiometry to artificially compensate for missing transient dynamics later in the discovery phase.

### 2.3.2. Isolation of Primary Monod Parameters ($\mu_{max}$ and $K_s$)

Under continuous, uninhibited operational baselines, the time derivatives of the bioreactor mass balance equate to zero, forcing the specific growth rate to exactly match the operating dilution rate ($\mu = D$). By evaluating the steady-state substrate concentrations ($S$) across a controlled range of baseline dilution rates, the primary Monod kinetics are entirely isolated from any secondary regulatory feedback.

$$\mu = \frac{\mu_{max} S}{K_s + S} = D \tag{9}$$

Through nonlinear regression of this steady-state baseline data, the maximum specific growth rate ($\mu_{max}$) and the substrate affinity constant ($K_s$) are explicitly calculated. Fixing $\mu_{max}$ and $K_s$ establishes a rigid primary growth manifold, ensuring the foundational carbon-uptake logic remains stable during subsequent dynamic analysis.

### 2.3.3. Isolation of the Gated-Hill Yield Parameters ($k_{base}$, $k_{max}$, $K, n$)

A similar steady-state decoupling is applied to the secondary metabolite yield ($Y_{P_s/X}$). The parameters governing the Gated-Hill architecture define the physical and regulatory boundaries of the metabolic shift. By systematically varying the steady-state dilution rate to walk the reactor through different stable substrate concentrations, the physical yield curve is explicitly mapped.

Because the thermodynamic gate ($K_{gate}$) is strictly constrained to a physical boundary limit, the remaining parameters—the basal production rate ($k_{base}$), the maximal derepressed rate ($k_{max}$), metabolic shift sharpness (n) and the macroscopic starvation threshold ($K$)—can be isolated and optimized purely from the steady-state yield variance.

By resolving and locking these baseline kinetic and yield parameters prior to dynamic analysis, the model discovery engine is strictly constrained. When the dynamic system data is transformed into the log domain to convert multiplicative kinetic models into additive linear features, the presence of these fixed anchors prevents the algorithm from minimizing errors through fractional weight distributions or physically impossible parameter combinations. With the fundamental stoichiometric mass balances and growth boundaries locked, the discovery algorithm is forced to allocate all transient state-space variance exclusively to identifying the unknown dynamic mechanisms, successfully isolating complex phenomena such as the product toxicity shutdown limits and the autocatalytic positive feedback loop.

## 3. THE AUTONOMOUS KINETIC MANIFOLD DISCOVERY FRAMEWORK

The Autonomous kinetic manifold discovery framework is a closed-loop diagnostic architecture designed to autonomously detect, isolate, and reconstruct missing physical dynamics within biological systems. Its elements and the information flow among them are illustrated in Figure 3. Step-by-step explanation of the workflow is given next.

The workflow begins in the Offline (Design of Experiment) DoE Framework box. This component designs experimental protocols for the process inputs which force the reactor to different regimes. This plan is created by the scientist with the help of the agent. The Scientist takes the offline plan and, through arrow (1), "Defines Objectives & Known Priors ($\mu_{max}$ , $K_s$). This sets the boundary conditions and starting points for the AI system (LLM Discovery Agent).

The LLM Discovery Agent functions as the core "Scientific Discovery Agent," and initiates the search for candidate mathematical models. It queries the Two-Level RAG Knowledge Base, which is structured with both multiplicative and additive models. After the query is complete, the LLM Discovery Agent must internally "process" the results. This is where the Agent's internal "JSON Physics Generator" selects the most promising model candidate from the RAG search, formats it as a mathematical description (specifically an ODE), and prepares it for simulation.

The Agent now transfers the formulated mathematical model to the simulation engine. Arrow (4), "JSON ODE Blueprint," represents the transfer of the ODE model structured in JSON format from the Agent down to the MATLAB Execution Bridge. This blueprint includes instructions like build ODEs. The MATLAB Execution Bridge acts as the laboratory emulator. It takes the "JSON ODE Blueprint" and combines it with physical system data. This data is the "Cleaned Trajectory Data (X, S, P, $\mu_{obs}$)" generated from the Physical & Data Environment, which processes raw CSTR plant dynamics through MHE state estimation. MATLAB performs a "Tournament Optimization," testing multiple versions of the model against the cleaned data to find the best fit.

The results of the optimization tournament are calculated as statistical criteria (e.g., Akaike Information Criterion, or AIC scores). Arrow (5), "AIC Scores & Tournament Results," sends these values back up to the LLM Discovery Agent for evaluation. The LLM Agent's "Degeneracy Analyzer" reviews the AIC scores. If competing models are mathematically "degenerate" (indistinguishable, providing identical fits), the "Active DoE Controller" decides to intervene. This triggers Step (6): the red dashed line labeled "Active DoE Pulse (Triggered on ΔAIC Degeneracy)" sends an instruction back to the physical CSTR Plant to perform an intervention, such as changing substrate concentration (S). These steps (indicated by 2-6) then form a recursive loop that continues until the model degeneracy is resolved.

Once the LLM Agent resolves model degeneracy and finds a unique, validated model that explains the multi-modal physics of the reactor, the process concludes. Arrow (7), "Validated Digital Twin," sends the final, confirmed model, presented as a "Final Reconstructed Manifold Dashboard," back up to the Scientist / Operator for review, approval, and final acceptance.

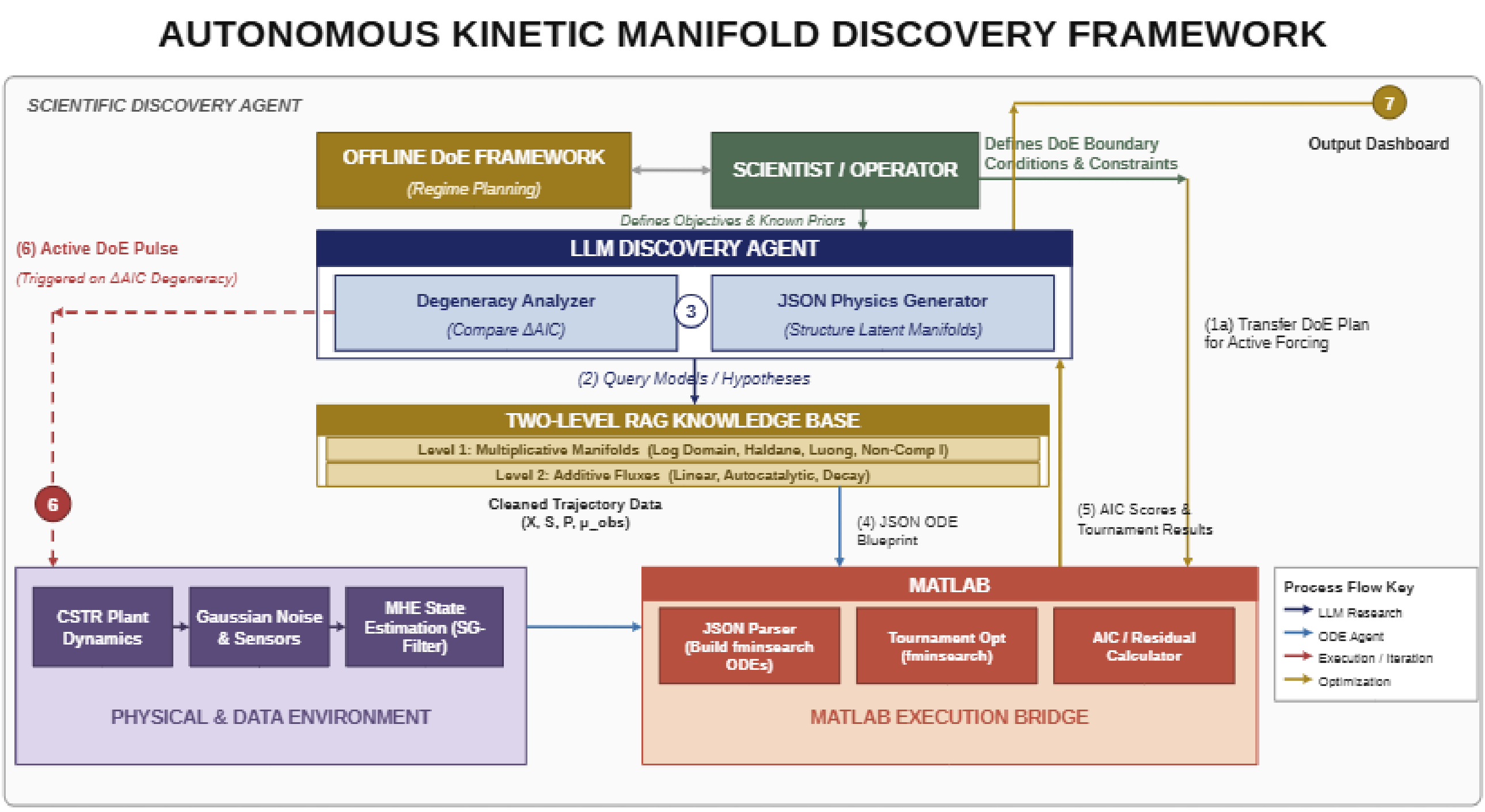

**Figure 3**. Autonomous Discovery Architecture.

Note that DoE is performed both offline (by the Scientist) and online (by the LLM). While a fully autonomous, LLM-driven operational sequence is theoretically appealing, deploying an LLM to manage routine regime transitions online introduces unnecessary computational overhead and tuning fragility. Therefore, this framework adopts a hybrid control strategy: classical, offline DoE is utilized to establish the baseline orthogonal regimes, reserving the LLM exclusively as an "active experimenter" to intervene only when the framework encounters mathematical degeneracy.

We now explain the procedural details used by different members of this discovery architecture.

### 3.1 Moving Horizon State Estimation (MHE)

To isolate the structural residuals necessary for kinetic discovery, the system must first extract the time-varying, true specific growth rate ($\mu_{obs}$) from the noisy reactor state trajectory. Because the governing biological kinetics are treated as an unknown latent manifold, we utilize a Moving Horizon Estimation (MHE) protocol to continuously track $\mu$ prior to any structural identification.

The state estimator operates on the discrete state history obtained from the CSTR and isolates the biomass trajectory $X_{raw}$. To mitigate the impact of Gaussian measurement noise ($\sigma = 0.05$) and prevent high-frequency noise from artifactually triggering the active discovery agent, the raw biomass signal is pre-processed using a Savitzky-Golay filter with a rolling window of 250 discrete data samples (corresponding to a physical temporal window of 125 hours, given the sampling interval of $\Delta t = 0.5$ h).

The core estimation is performed across a rolling computational window of $N_{horizon} = 10$ time steps. At each discrete time $k$, the estimator utilizes Sequential Quadratic Programming (SQP) via fmincon to find a local constant $\mu_k$ that best reconstructs the smoothed biomass trajectory. The optimization minimizes the sum of squared errors between the filtered data and a local Euler integration of the biomass mass balance, governed by the objective function:

$$J(\mu_k) = \sum_{i=k-N_{horizon}+1}^{k} \left(X_{sg}(t_i) - X_{sim}(\mu_k, t_i)\right)^2 + \lambda(\mu_k - \mu_{k-1})^2 \quad (10)$$

To enforce physical continuity and prevent the estimator from tracking non-physical, high-frequency oscillations, the objective function incorporates a Tikhonov regularization penalty. A tuning parameter of $\lambda = 0.10$ heavily penalizes discontinuous jumps between the current estimate $\mu_k$ and the previous estimate $\mu_{k-1}$. The output of this MHE protocol is a continuous, mathematically robust state trajectory for the total specific growth rate, $\mu_{obs}$.

### 3.2 Logarithmic Transformation

For the biological CSTR under consideration, the specific growth rate is expressed as:

$$\mu = \underbrace{\left(\frac{\mu_{max}S}{K_s+S} \cdot \mu_H(S) \cdot \mu_L(P)\right)}_{\text{Multiplicative Constraint Path}} + \underbrace{\mu_{auto}}_{\text{Additive Rescue Path}} \tag{11}$$

We assume that the foundational Monod kinetic parameters—specifically the maximum specific growth rate ($\mu_{max}$) and the half-saturation constant ($K_s$)—are known and fixed from available kinetic data collected during the reactor operation in the low-substrate Monod regime. Consequently, the primary challenge is to autonomously discover the physics governing the deviations beyond standard Monod kinetics.

The specific growth rate expression (Eq. 11) is composed of a multiplicative and additive term. While logarithmic transformation elegantly decomposes multiplicative kinetic structures, it cannot directly resolve additive kinetic components. Mathematically, the transformation $\ln(\mu_{base} \cdot \mu_H \cdot \mu_L + \mu_{auto})$ does not yield linearly separable terms since $\ln(A \cdot B + C) \neq \ln(A) + \ln(B) + \ln(C)$.

To overcome this limitation, the discovery architecture leverages the chronological decoupling inherent to bioprocess dynamics. Secondary additive pathways, such as autocatalytic maintenance, typically remain dormant or negligible during early exponential growth and primary steady-state operations ($\mu_{auto} \approx 0$). By sequentially partitioning the discovery tournament, the agent first operates within a temporal window where the kinetic manifold is strictly multiplicative. During this phase, the MHE estimate represents the isolated multiplicative growth response ($\mu_{obs,mult}$). This allows the system to decompose the multiplicative log residual to the individual structural inhibitory residuals:

$$E_{res} \equiv \ln(\mu_{obs,mult}) - \ln(\mu_{base}) = \ln(\mu_H) + \ln(\mu_L) \quad (12)$$

After the individual terms are discovered and reconstructed into a cumulative multiplicative estimate ($\mu_{recon_mult}$), the final additive residual is isolated via straightforward subtraction against the total observed specific growth rate ($\mu_{obs}$), enabling the subsequent discovery of the delayed autocatalytic parameters without mathematical interference:

$$\epsilon_{add} = \mu_{obs} - \mu_{recon_mult} \quad (13)$$

### 3.3 The Principle of Orthogonal Additive Regimes

A fundamental challenge in bioprocess system identification is the severe mathematical coupling between competing kinetic mechanisms. In a standard batch or unoptimized continuous trajectory, substrate inhibition, product toxicity, and secondary metabolic fluxes occur simultaneously, creating a highly entangled kinetic landscape. To overcome this, our discovery architecture relies on a Design of Experiments (DoE) framework that actively forces the reactor into *orthogonal additive* regimes. By exploiting the logarithmic projection of the multiplicative growth manifold, inhibitory mechanisms are mathematically transformed into a linear, additive superposition:

$$\ln(\mu_{obs,mult}) = \ln(\mu_{base}) + \ln(\mu_H) + \ln(\mu_L) \quad (14)$$

For these terms to be uniquely identifiable, the DoE protocol utilizes dynamic closed-loop control to temporarily suppress specific state variables (driving their corresponding logarithmic residuals to zero) while intentionally exciting others. This ensures that during any given discovery phase, the isolated structural residual is orthogonally decoupled from the rest of the manifold, guaranteeing that the AI agent optimizes against a single, unobstructed physical mechanism.

### 3.4 Controlled Protocols for Regime Isolation

To systematically reconstruct the complete kinetic state-space, the Design of Experiments (DoE) applies a series of orthogonal step-setpoint changes to transition the reactor through seven distinct operating regimes. By forcing the continuous stirred-tank reactor (CSTR) into targeted steady-state spaces, structural model mismatches are converted into isolated, sustained flatline residual signatures. Once these baselines are established, highly constrained, short-duration transient pulses are superimposed strictly during the final generation step to excite localized kinetic sensitivities.

To maintain optimal operational regimes within the continuous stirred-tank reactor (CSTR) while simultaneously performing Moving Horizon Estimation (MHE), a robust and decoupled control strategy is required. Classical approaches to CSTR control often utilize model-based control, relying on the estimated specific growth rate ($\mu_{est}$) to predict biological substrate consumption and adjust feed rates accordingly. However, incorporating real-time $\mu_{est}$ directly into the control law introduces an algebraic loop between the estimator and the controller. This couples the systems, amplifying sensor noise, destabilizing the numerical integration, and causing significant high-frequency chattering in the manipulated variables.
To resolve this, we propose a purely error-driven Proportional (P) control scheme. There is no need for integral action since off-set free tracking is not necessary as long as the states enter the particular regimes of interest. The control laws for the dilution rate ($D$), inlet substrate concentration ($S_{in}$), and inlet product concentration ($P_{in}$) are defined as follows:

$$D = D_{bias} + K_{p,X}(X_{sp} - X_{filt}) \tag{15}$$

$$S_{in} = S + K_{p,S}(S_{sp} - S_{filt}) \tag{16}$$

$$P_{in} = P + K_{p,P}(P_{sp} - P_{filt}) \tag{17}$$

States that are measured with noise are filtered by an exponential moving average filter and the filtered values ($X_{filt}, S_{filt}, P_{filt}$) are fed to the controllers. Unlike standard industrial controllers that rely on a fixed nominal baseline (e.g., $S_{in,nominal}, P_{in,nominal}$), the feed pump controllers for $S$ and $P$ dynamically utilize the current filtered state ($S_{filt}$ *and* $P_{filt}$) as their baseline bias. The physical justification for this state-dependent baseline becomes evident when substituting the control law into the standard CSTR substrate mass balance:

$$\frac{dS}{dt} = D(S_{in} - S) - \frac{\mu X}{Y_{X/S}} \tag{18}$$

Substituting the proposed $S_{in}$controller into the mass balance yields:

$$\frac{dS}{dt} = D\left(S + K_{p,S}\left(S_{sp} - S\right) - S\right) - \frac{\mu X}{Y_{X/S}} \tag{19}$$

The $+S$ and $-S$ terms perfectly annihilate each other, which successfully neutralizes the reactor's natural washout penalty ($-DS$). In a standard CSTR, higher substrate concentrations are physically flushed out at a faster absolute rate, forcing a traditional controller to constantly fight the baseline hydrodynamics of the tank. By eliminating this specific nonlinear washout coupling, the controller output acts as a more direct, stable steering mechanism, heavily dampening transient oscillations. The sequence of regimes is formalized in Figure 4 (a,b,c).

Initially (0 to 50 hours) the system operates in the baseline monod conditions achieved by Setpoints $X_{sp} = 0.05, S_{sp} = 20.0, P_{sp} = 0.0$. This establishes initial hydraulic and metabolic equilibrium. This is followed by 3 distinct regimes:

- **Regime I: Substrate Inhibition Isolation (50 to 150 hours):** Setpoints $X_{sp} = 0.05, S_{sp}\ is\ increased\ to\ 30.0\ at\ t = 50\ and\ decreased\ back\ to\ 20\ at\ t = 150, P_{sp} = 0.0$. Elevated substrate levels map the Haldane structural threshold.
- **Regime II: Product Toxicity Isolation (250 to 500 hours):** Setpoints $X_{sp} = 0.05, S_{sp} = 20.0, P_{sp} = increased\ from\ 0\ to\ 40.0$ at t=250 and decreased to 0 at t=320. Maps multiplicative metabolic suppression caused by products.
- **Regime III: Metabolic Starvation (500 to 800 hours):** Setpoints $X_{sp} = 0.05, S_{sp} is\ dropped\ from\ 20\ to\ 2\ to\ enter\ metabolic\ starvation\ where\ netabolite\ Ps\ is\ produced.$

  At t=800 the reactor is returned to its base line.

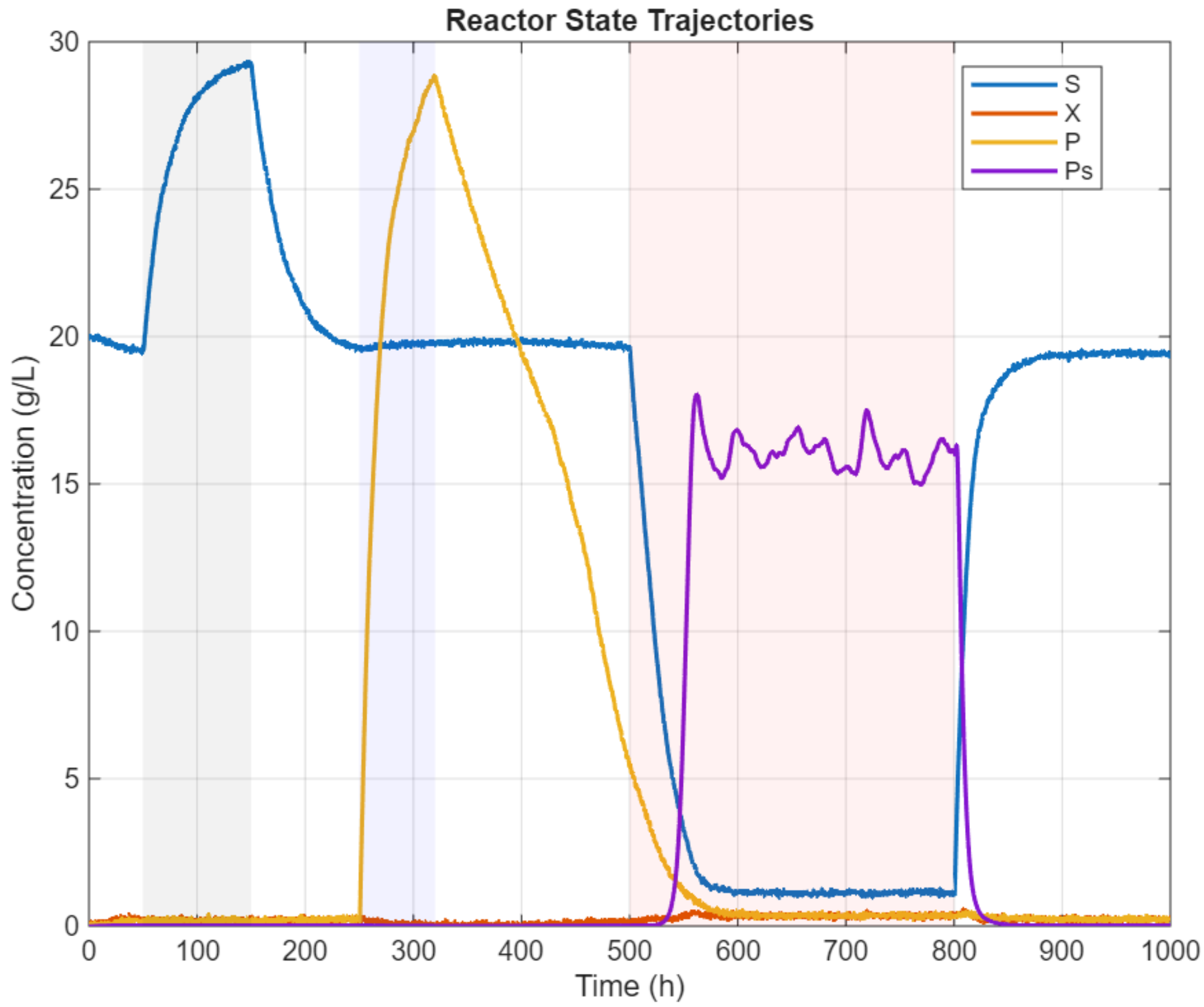

**Figure 4.a**. The operating regimes used for discovery.

Note that in Regime III when substrate concentration is decreased to enter starvation, metabolite production responds to this set-point change in a somewhat oscillatory fashion. Since the objective is to be able to discover the metabolite, no attempt is made in this regime to control it directly. By intentionally leaving $P_s$ uncontrolled, we allow it to act as a free-floating, dynamic state variable, generating the rich, unforced variance that the Autonomous Agent needs in Phase 3 to correctly map the hidden autocatalytic ($k_{auto}$) parameters. The transition from one regime to another is achieved by the controllers. Manipulated variables are shown in Figure 4b and c.

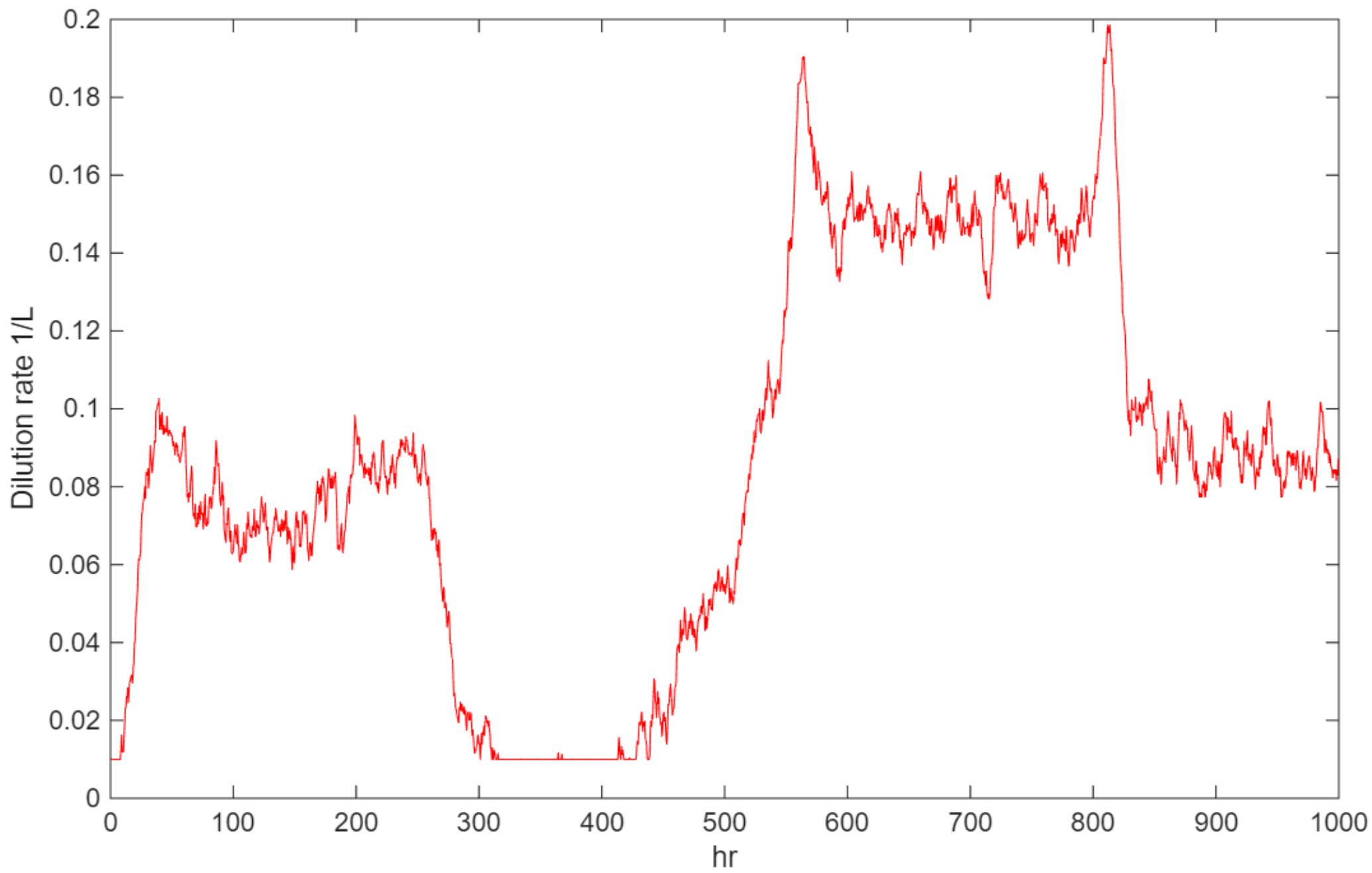

**Figure 4.b**. Dilution rate during regime transitions.

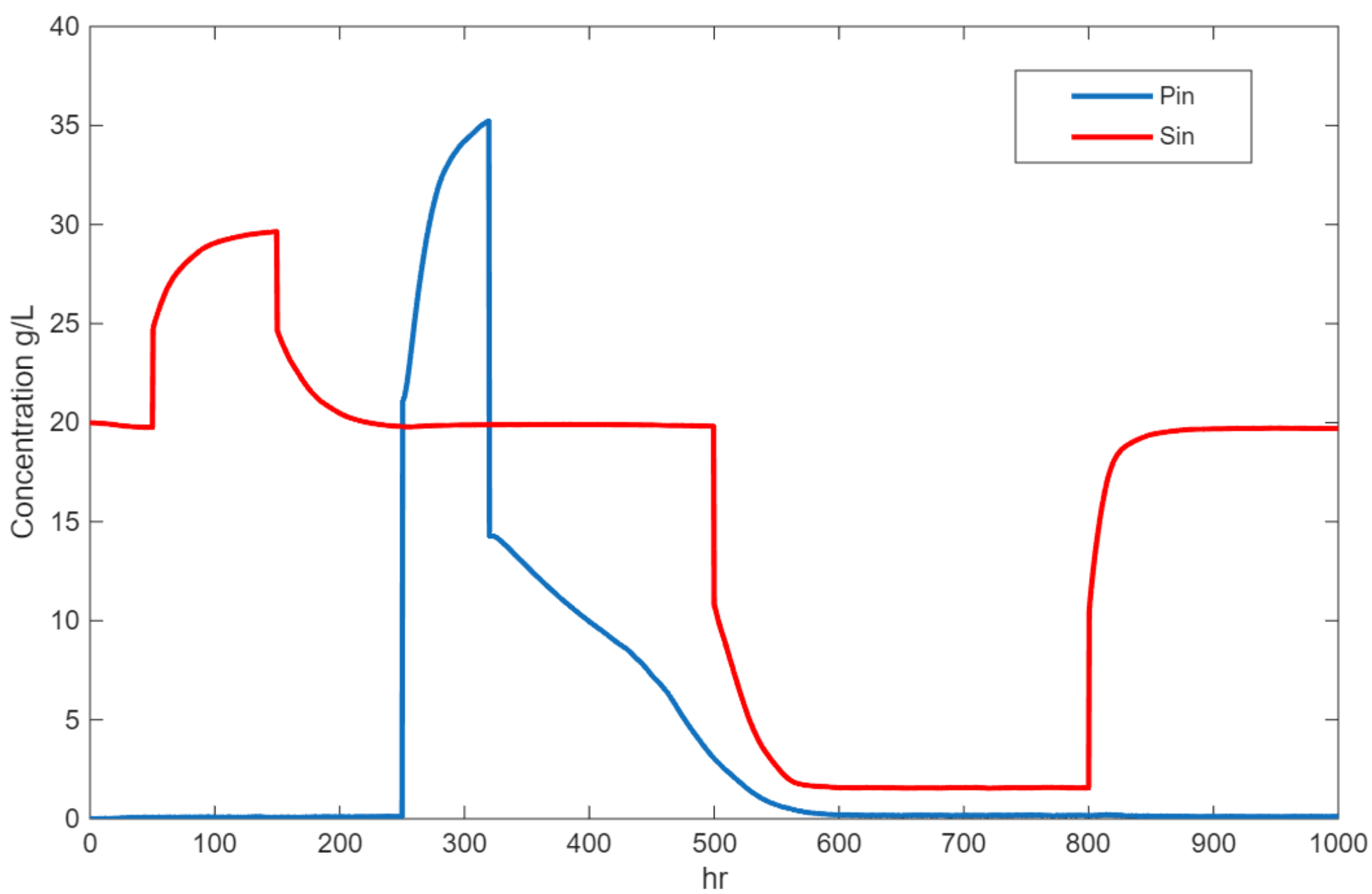

**Figure 4.c.** Feed concentrations during regime changes.

### 3.5 The Discovery Engine

The discovery engine acknowledges the inherent plant-model mismatch by treating structural complexity as an unknown variable to be recovered through autonomous observation and diagnostic perturbation. The Discovery Engine processes the data through two distinct Discovery Stages, each tied to specific operational regimes of the CSTR.

**Stage 1: Mapping the Primary Manifold.** During the initial operational window ($t < 500$ h), the agent actively probes the reactor via controlled substrate and product shifts to isolate the underlying multiplicative kinetic terms. Operating strictly within the logarithmic domain, it fits the Haldane inhibition model to the substrate ($S$) and the Luong toxicity model to the product ($P$). Upon completing these fits, the agent reconstructs the predicted multiplicative growth rate ($\mu_{inhibitory}$) and subtracts it from the true biological growth rate identified by the Moving Horizon Estimator ($\mu_{tracked}$). The algorithm then evaluates the Mean Residual Variance (MRV) of this unexplained error. If the MRV is sufficiently low—approaching the sensor noise floor—and the residuals exhibit random, zero-mean behavior, the agent mathematically confirms that the primary $(S, P)$ state-space is exhaustively mapped. At this point, no further structural signal can be extracted from the primary variables.

**Stage 2: Transition to Additive Search.** As the reactor is deliberately forced into the starvation regime ($t \geq 500$ h), a sudden and substantial spike in the residual error emerges. Because the agent mathematically verified in Stage 1 that the primary states ($S$ and $P$) cannot account for this new variance, it formally terminates the multiplicative log-domain search. The pipeline then drops out of the logarithmic space and transitions to evaluating secondary, unmodeled sensor data as additive terms to explain this newly unmasked error block.

**RAG Library**

By integrating a two-level Retrieval-Augmented Generation (RAG) library with an autonomous Large Language Model (LLM) agent, the discovery engine iteratively evolves from a simplified Monod baseline to a full-rank, five-dimensional state-space representation. Because biological growth features both inhibitory (multiplicative) and supplementary (additive) mechanisms, the discovery architecture uses a Two-Level Dual-Library System (see Table 1) structured as follows:

- **Level 1 Library (Multiplicative / Log Domain):** Contains only the gating and inhibitory functions (Shelves 1, 2, and 3). These are log-transformed kinetic basis functions designed to scale the primary baseline growth.
- **Level 2 Library (Additive / Linear Domain):** Contains only the auxiliary metabolic fluxes (Shelf 4), such as metabolic sources (e.g., precursor-driven autocatalysis) and metabolic sinks (e.g., endogenous decay, maintenance energy, and non-growth-associated Luedeking-Piret production). These are strictly evaluated in the linear domain after Level 1 mapping is complete.

**The Additive Illusion**

If a true multiplicative phenomenon governs the system, but the corresponding structure is absent from the Level 1 Library, the agent may erroneously categorize the unmodeled variance as an additive error upon transitioning out of the logarithmic domain. While the agent cannot intrinsically detect this misclassification, a human supervisor can identify the structural illusion through specific mathematical heuristics.

A genuine additive term (e.g., background secondary metabolite production) should operate independently of the primary baseline growth rate. Conversely, the residual error generated by a missing multiplicative term remains permanently scaled by the primary baseline growth, $\mu_{base}$. Therefore, if the unexplained linear variance perfectly tracks the transient oscillations of the observed growth rate ($\mu_{obs}$) or the primary substrate ($S$), the missing physics are almost multiplicative.

Furthermore, a falsely identified additive term will be severely overfit to the specific operational window of its discovery. If the baseline operating conditions are altered—such as doubling the substrate setpoint ($S_{sp}$) or radically shifting the dilution rate ($D$)—the primary growth will scale, but the artificial additive term will fail to scale proportionally, causing an immediate collapse in the model's predictive validity. Finally, if the linear error is plotted against a suspected secondary additive variable and the resulting line of best fit does not pass near the origin $(0,0)$—indicating substantial residual error even when the secondary state is zero—it serves as a definitive signal of structural mismatch.

Ultimately, an autonomous agent is strictly bound by the mathematical vocabulary it is provided. If the true underlying physics do not exist within the library, the algorithm is mathematically forced to construct a structural fiction from the pieces it has left.

**Table 1:** The Physical Retrieval-Augmented Generation (RAG) Library for CSTR Kinetic Manifolds

| RAG Shelf / Phase | Model Name | Standard Physical Equation | AI Solver Target String | Citation | Typical Substrate / Microorganism / Product |
|---|---|---|---|---|---|
| **Shelf 1: Substrate Inhibition (Multiplicative Log Domain)** | Haldane | $\frac{1}{1+S^2/(K_i \cdot (K_s+S))}$ | $-\ln(1+S^2/(K_i \cdot (K_s+S)))$ | Andrews (1968) | Phenol / Pseudomonas putida / Biomass |
| | Luong Substrate | $(1-S/S_{max})^n$ | $n \cdot \ln(1-S/S_{max})$ | Luong (1987) | Butanol / Clostridium acetobutylicum / Solvents |
| | Aiba (Substrate) | $\exp(-S/K_i)$ | $-S/K_i$ | Aiba et al. (1968) | Glucose / Saccharomyces cerevisiae / Ethanol |
| | Webb (Survival) | $\frac{1+S/K_w}{1+S/K_i}$ | $\ln(1+S/K_w)-\ln(1+S/K_i)$ | Webb (1963) | Enzyme kinetics (Toxic substrate with survival) |
| | Yano & Koga | $\frac{1}{1+S/K_1+S^2/K_2}$ | $-\ln(1+S/K_1+S^2/K_2)$ | Yano & Koga (1969) | Toluene / Pseudomonas sp. / Biomass |
| **Shelf 2: Product Toxicity (Multiplicative Log Domain)** | Luong Toxicity | $(1-P/P_m)^n$ | $n \cdot \ln(1-P/P_m)$ | Luong (1985) | Glucose / S. cerevisiae / Ethanol |
| | Aiba (Product) | $\exp(-kP)$ | $-k \cdot P$ | Aiba et al. (1968) | Glucose / S. cerevisiae / Ethanol |
| | Jerusalimsky | $\frac{1}{1+P/K_p}$ | $-\ln(1+P/K_p)$ | Jerusalimsky (1967) | Glucose / Propionibacteriu |

| RAG Shelf / Phase | Model Name | Standard Physical Equation | AI Solver Target String | Citation | Typical Substrate / Microorganism / Product |
|---|---|---|---|---|---|
| | | | | | m / Propionic Acid |
| | Hinshelwood | $(1 - kP)$ | $\ln(1 - kP)$ | Hinshelwood (1946) | Generic bacterial cultures |
| **Shelf 3: Causal Novelty (Multiplicative Log Domain)** | Non-Comp. Byproduct (Latent Biotoxin, B) | $\frac{1}{1 + B/K_{IB}}$ | $-\ln(1 + B/K_{IB})$ | General Enzyme Kinetics | Glucose / Escherichia coli / Acetate (Byproduct) |
| | Competitive Byproduct | $\frac{1}{1 + (B/K_{IB}) \cdot (K_s/(K_s + S))}$ | $-\ln(1 + (B/K_{IB}) \cdot (K_s /(K_s + S)))$ | General Enzyme Kinetics | Xylose / S. cerevisiae / Xylitol |
| | Unmodeled Toxic Substrate | $\exp(-k \cdot S_{unmodeled})$ | $-k \cdot S_{unmodeled}$ | Environmental Engineering | Mixed VOCs / Biofilter Consortia / Biomass |
| | Substrate/Product Synergy | $\exp(-k \cdot S \cdot P)$ | $-k \cdot S \cdot P$ | Levenspiel (1980) | Complex high-gravity fermentations |
| **Shelf 4: Additive Fluxes (Additive Linear Domain)** | Null Hypothesis (No Secondary Effect) | 0 | 0 | Base Statistical Control | Pure Gaussian Noise / Null Error Baseline |
| | Linear Product Boost (Non-Saturating Autocatalysis) | $k_{linear} \cdot P_s$ | $k_{linear} \cdot P_s$ | Simplified QS Models | AI-2 / Vibrio harveyi / Target Gene Expression |
| | Saturating Product Boost (Monod Autocatalysis) | $k_{sat} \cdot \left(\frac{P_s}{K_{ps} + P_s}\right)$ | $k_{sat} \cdot \left(\frac{P_s}{K_{ps} + P_s}\right)$ | Monod (1949) adapted for Quorum Sensing | AHL (Autoinducer) / Vibrio fischeri / Bioluminescence |
| | Hill/Cooperative Signaling | $k_{auto} \cdot \frac{P_s^n}{K_p^n + P_s^n}$ | $k_{auto} \cdot \frac{P_s^n}{K_p^n + P_s^n}$ | Hill (1910) | Auto-inducer Peptide (AIP) / Staphylococcus aureus / Virulence |
| | Constant Additive Leakage | $k_{auto}$ | $k_{auto}$ | Synthetic Biology Baselines | Engineered Promoters / E. coli / |

| RAG Shelf / Phase | Model Name | Standard Physical Equation | AI Solver Target String | Citation | Typical Substrate / Microorganism / Product |
|---|---|---|---|---|---|
| | | | | | Recombinant Protein |
| | Endogenous Decay | $-K_d$ | $-\lvert K_d \rvert$ | Herbert (1958) | Universal Biomass Turnover / Maintenance |
| | Density Lysis | $-K_{lysis} \cdot X$ | $-\lvert K_{lysis} \rvert \cdot X$ | | |

### 3.6 Structural Identification

Once the structural residual ($E_{res}$) is isolated following the Moving Horizon Estimation (MHE) and logarithmic transformation protocols, the LLM agent first queries the Level 1 RAG Library. As detailed in Table 1, this primary library exclusively houses candidate functional forms for multiplicative kinetic modifiers, encompassing mechanisms such as Luong toxicity, Jerusalimsky inhibition, and various competitive or non-competitive interactions. To reconstruct the missing physics, the discovery engine executes a computational tournament, mathematically evaluating these competing hypotheses against the isolated residual to seamlessly identify and align the true structural match.

For each candidate kinetic basis function $f_j(\theta, \mathbf{x})$ retrieved from the RAG library, the agent formulates a Sum of Squared Errors (SSE) objective function $J(\theta)$. This function quantifies the deviation between the candidate model and the isolated structural residual $E_{res}$ across $N$ discrete time steps within the targeted operational phase:

$$J(\theta) = \sum_{k=1}^{N} \left( E_{res}(t_k) - f_j\big(\theta, \mathbf{x}(t_k)\big) \right)^2 \quad (20)$$

where $\theta$ represents the vector of unknown kinetic parameters (e.g., $[K_{ib}]$ or $[P_m, n]$) for the candidate model, and $\mathbf{x}(t_k) = [X(t_k), S(t_k), P(t_k), Ps(t_k)]^T$ represents the corresponding physical state vector of the reactor.

It is important to note that during Stage 1, the agent's active search space is strictly restricted to the primary dimensions $[S, X, P]$. While the autocatalytic precursor $P_s$ is continuously produced and present within the physical reactor, the agent deliberately excludes this sensor data from its observable state-vector. $P_s$ is integrated into the agent's active system identification only when the algorithm explicitly hypothesizes the presence of secondary metabolites in Stage 2 to explain the additive residual error.

By minimizing this objective function, the optimization algorithm actively searches for the specific physical constants ($\theta_j$) that best align the empirical data trajectory with the fundamental geometric shape of the hypothesized basis function:

$$\hat{\theta}_j = \arg\min_{\theta} J(\theta) \quad (21)$$

To solve this unconstrained optimization problem, the discovery agent relies strictly on a gradient-free Nelder-Mead simplex search method. Pivoting away from gradient-reliant algorithms is a deliberate architectural constraint, preventing the solver from becoming trapped in the highly non-convex, localized topology inherent to complex, overlapping biological kinetic equations.

To prevent mathematical overfitting—where a highly complex, non-physical polynomial might artificially outscore a true biological mechanism—candidates are strictly evaluated using the Akaike Information Criterion (AIC).

$$\text{AIC} = N \cdot \ln\left(\frac{\text{SSE}}{N}\right) + 2k \quad (22)$$

where $N$is the number of observations, SSE is the sum of squared errors, and $k$is the number of estimated parameters. The AIC penalizes unnecessary complexity, ensuring the agent favors elegant, physically grounded mechanisms.

Once the LLM has successfully reconstructed the complete multiplicative backbone ($\mu_{mult}$), it shifts to the Level 2 RAG Library to discover any remaining auxiliary, growth-linked metabolic pathways (e.g., the autocatalytic maintenance loop driven by precursor $Ps$). Since these mechanisms act as parallel energy sources, they are mathematically additive. The logarithmic decomposition is no longer valid. Instead, the agent executes a direct linear projection to isolate the final temporal residual:

$$\epsilon_{add} = \mu_{obs} - \mu_{recon_mult} \tag{23}$$

The agent then queries the Level 2 library, which houses strictly additive kinetic structures. Fitting candidate models to $\epsilon_{add}$ allows the precise identification of the autocatalytic rate constants without any mathematical interference from the previously resolved inhibitory dynamics. The LLM is programmed to monitor the informational divergence between top candidates. If the difference between competing hypotheses falls below a critical threshold ($\Delta\text{AIC} < 2.0$), the agent flags an "Identifiability Trap," recognizing that the data lacks the dynamic variance necessary to differentiate physics. To break this degeneracy, the LLM seamlessly transitions from a passive observer to an Active Experimenter. It autonomously commands an online Design of Experiments (DoE) intervention. By injecting a high-magnitude targeted substrate impulse (+15 g/L at t = 360 h), the agent dynamically decorrelates the competing model structures. The RAG tournament is then re-executed on the newly generated transient data, decisively resolving the structural ambiguity and locking in the true parameter.

## 4. RESULTS AND DISCUSSION

The autonomous discovery pipeline was executed against the hidden biological plant simulation. The digital exhaust from the agent's execution—detailing model selection, statistical confidence metrics (AIC), and autonomous intervention logic is presented below.

Console execution log of the autonomous LLM RAG Discovery Agent, detailing the iterative resolution of the primary multiplicative manifold (Level 1) and the active unmasking of the secondary additive flux (Level 2) is shown in Figure 5.

| Discovery Stage | Process Description & Agent Logic | Tested Hypotheses (Models) | Statistical Metrics (AIC/MRV) | Identified Parameters & Outcome |
|---|---|---|---|---|
| **1. System Initialization** | Agent Status: Active. Live Plant Simulation initialization and Sensor Data Logging complete. | | | |
| **Stage 1: Substrate Inhibition Tournament (Log Domain)** | Primary mechanism identification (substrate-driven inhibition).<br>Log Domain isolation. | Haldane | AIC −1034.17, MRV: 0.018874 | **Winner: Haldane Inhibition**<br>Identified Parameter: $K_i$ = 5.0940 |
| | | Aiba (Substrate) | AIC −990.80, MRV: 0.022286 | |
| | | Webb (Survival) | AIC −1024.65, MRV: 0.019425 | |
| | | Yano & Koga | AIC −1032.17, MRV: 0.018874 | |
| **Stage 2: Toxicity Screening Tournament** | Identification of high-product toxicity manifold.<br>Comparison of non-linear models. | Luong_Tox | AIC −144.65, MRV: 0.249262 | **Winner: Luong Toxicity Manifold**<br>Identified Parameters:<br>$P_m = 29.7206$ g/L, $n = 1.3191$<br>$n = 1.3191$ |
| | | Jerusalimsky | AIC −116.21, MRV: 0.331282 | |
| | | Aiba_Tox | AIC −132.85, MRV: 0.283583 | |
| | | Hinshelwood | AIC −143.09, MRV: 0.257691 | |
| **Level 1 Completion Verification** | Evaluating E_global across primary operational window (t < 500 h).<br>Evaluating residual error (E_global). | Mean Residual Variance = 0.0011 | | Residuals exhibit zero-mean behavior. Primary multiplicative manifold (S, P) mapped.<br>**AGENT ACTION:** Terminating Level 1. Initiating Secondary Sensor Evaluation. |
| **Stage 3: Secondary Metabolite Tournament (Additive Domain) (t = 500–750 h)** | 1. Primary state-space [S, P] mapped and exhausted. Residual error detected.<br>2. Evaluating measured side product (Ps) against unexplained residual error…<br>3. Identifiability Trap Detected (Linear vs. Saturating collinearity at low Ps).<br>4. Autonomous Active DoE Pulse Triggered (S_sp = 15.0 at t=700). (Preserve DoE detail format).<br>5. LLM dynamically inserts dPs/dt ODE into mathematical library.<br>6. Re-executing RAG Tournament on transient DoE data. | No_Ps_Effect | AIC −2282.56, MRV: 0.010504 | **Winner: Saturating Secondary Metabolite Feedback (Ps).**<br><br>Identified Parameters:<br>$k_{auto}$ = 0.1596<br>$K_p$ = 6.5502<br>**SUCCESS:** Degeneracy resolved. True structural manifold locked. |
| | | Linear_Ps_Boost | AIC −3934.54, MRV: 0.000387 | |
| | | Saturating_Ps | AIC −4065.76, MRV: 0.000297 | |

**Figure 5.** The Agent Log (*Summary of the Autonomous Kinetic Model Discovery Tournament and Structural Identification Results, multi-stage autonomous kinetic discovery process, structural identification logic, and final model structural locking results. Stage 4 highlights the agent's ability to detect identifiability traps and autonomously trigger active DoE to resolve model collinearity and insert missing physical ODEs.)*

### <u>4.1: Stage 1- Substrate Inhibition Tournament (Log Domain)</u>

During the initial 500-hour operational window (Regimes 1-3 of Figure 5), the true plant simulation was engineered using pure, classic Haldane substrate kinetics. The agent operated strictly within the log-transformed domain to isolate multiplicative rate phenomena. Analyzing the performance of the four candidate hypotheses from the log in Figure 5 reveals why the Haldane model emerged as the clear winner:

1. **The Winner: Haldane (Perfect Plant Alignment)**
   - *Physical Mechanism*: The Haldane model assumes that at high concentrations, substrate molecules bind to a second, allosteric site on the enzyme, rendering it inactive. This is modeled mathematically using a squared inhibition term in the denominator: $\mu = \mu_{max} \frac{S}{K_s+S}\left(\frac{1}{1+S^2/K_i(K_s+S)}\right)$.
   - *Why it Won:* Since the underlying plant data was generated using this exact mathematical structure, the Haldane model represents the *true physics* of the system. The data naturally conformed to its bell-shaped curve without requiring extra correction terms. Using the unconstrained simplex search algorithm, the agent accurately estimated the primary inhibition constant ($K_i$) to be 5.0940, which is virtually identical to the true plant parameter of 5.0.
2. **The Near-Miss: Yano & Koga (Redundant Complexity)**
   - *Physical Mechanism:* Yano & Koga is a more complex, higher-order inhibition model (often $1 + S^2/K_1(K_s + S) + S^3/K_2(K_s + S)$). It physically accounts for scenarios with multi-site binding or extreme, catastrophic toxicity that causes a sharper "cliff-edge" drop in growth rate at very high substrate loads.
   - *Why it Lost:* In this specific scenario, the Yano & Koga model achieved the exact same data fit quality as Haldane (MRV = $0.018874$). However, it carries two parameters ($K_1, K_2$) instead of one ($K_i$). To ensure parsimony, the Akaike Information Criterion (AIC) adds a strict $+2$penalty for every extra parameter. Because Yano & Koga did not provide a better error reduction to justify its added complexity, its AIC score was penalized by 2 points (AIC $= -1032.17$) compared to Haldane (AIC $= -1034.17$), causing it to lose the tournament.
3. **The Over-Constrained Webb (Physical Implausibility)**
   - *Physical Mechanism*: The Webb model introduces a $\beta S$ term in the numerator to account for microbial "survival" or metabolic adaptation, suggesting the specific growth rate plateaus at a minimal level at very high substrate concentrations rather than crashing to zero.

- *Why it Lost:* The agent enforced a strict lower bound on the survival parameter $\beta$ to prevent non-physical behavior. While this allowed it to capture some data features, forcing the model to maintain this artificial survival floor caused its overall error variance to rise (MRV $= 0.019425$). This increased error variance, combined with the complexity penalty for two parameters, resulted in an AIC score of $-1024.65$, significantly worse than pure Haldane.

4. **The Structurally Inflexible Model: Aiba (Poor Fit)**
   - *Physical Mechanism:* The Aiba model uses an exponential inhibition form ($\mu = \mu_{max} \frac{S}{K_s+S} e^{-S/K_i}$). Physically, this implies a different mechanism of interaction, typically not associated with enzyme-substrate complexation dynamics.
   - *Why it Lost:* The exponential decay structure of Aiba lacks the rational polynomial denominator flexibility of Haldane. It could not adequately capture the specific shape of the substrate inhibition curve in this dataset, resulting in a substantially higher residual error (MRV = $0.022286$) and the worst AIC score in the tournament.

## 4.2: Stage 2- Toxicity Screening Tournament (Product Domain)

Here is the detailed structural and statistical rationale for the outcomes of the secondary screening tournament:

1. **The Winner: Luong (Capturing Non-Linear Toxicity)**

- *Physical Mechanism:* The Luong model establishes a hard biological limit by defining a maximum product concentration ($P_m$) at which cellular growth completely ceases. More importantly, it introduces a dimensionless exponent ($n$) that dictates the geometric curvature of the growth decline, allowing it to model accelerating or decelerating toxicity rates as the product accumulates.

- *Why it Won*: Unlike Stage 1—where extra parameters were rigorously punished by the AIC because the true data lacked complexity—Stage 2's underlying plant dataset features a genuinely non-linear toxicity exponent ($n = 1.42$). By utilizing its two-parameter flexibility ($P_m$ and $n$), the Luong model allowed the optimizer to find $n = 1.3191$. This structural alignment drastically cut down the mean residual variance (MRV: 0.249262). Even though the AIC formula applied a strict $+2$ penalty for carrying the second parameter, the massive drop in error variance successfully compensated for the penalty, allowing Luong to edge out its closest competitor.

2. **The Close Second: Hinshelwood (Rigid Linearity)**

- *Physical Mechanism*: The Hinshelwood model operates as a mathematically restricted version of Luong where the exponent is permanently fixed to exactly one ($n = 1$). This forces a strictly linear decline in the specific growth rate as product concentration increases.

- *Why it Lost*: Because it only requires a single parameter ($P_m$), Hinshelwood benefited from a lighter AIC penalty. However, its rigid linear assumption structurally failed to map the actual curvature of the plant's toxicity data. The increased error from this structural mismatch barely outweighed its parsimony advantage, causing it to lose the tournament by a narrow margin (AIC: -143.09 vs -144.65).

3. **The Structural Failures: Jerusalimsky and Aiba (Asymptotic Decay)**

- *Physical Mechanism*: Both the Jerusalimsky (rational saturation) and Aiba_Tox (exponential decay) models assume that product toxicity approaches zero asymptotically. Mathematically, they dictate that the culture will continue to maintain a fractional, non-zero growth rate infinitely, no matter how high the product concentration climbs.

- *Why they Lost:* In this specific bioprocess environment, the toxicity acts as a hard, absolute boundary—a cliff-edge where biological activity shuts down completely. Because Jerusalimsky and Aiba lack a mathematical mechanism to reach a true zero-growth state at a finite concentration, they consistently over-predicted growth at high product levels. This structural failure to capture the sharp asymptotic drop-off led to highly inflated residual errors (MRV: 0.331282 and 0.283583, respectively) and the poorest AIC standings in the tournament.

## 4.3 Stage 3 - Secondary Metabolite Tournament (Additive Domain)

Here is the detailed structural rationale for the final stage of the discovery pipeline:

**1. Secondary Sensor Evaluation and the Identifiability Trap**

In modern Process Analytical Technology (PAT) environments, tracking secondary metabolites like $P_s$ is standard practice, though their kinetic impact is often unknown to the baseline model. When the agent detected a structured, unexplained residual error remaining after mapping the primary $[S, P]$ state-space, it sequentially evaluated available secondary sensor streams to resolve the discrepancy. By utilizing the measured $P_s$ data as the input for a newly generated hypothesis library, the agent tested whether additive models could minimize the residual variance.

However, evaluating this measured secondary side product ($P_s$) against the residual error at low or steady-state concentrations revealed a severe identifiability trap. Under these baseline conditions, the mathematical signature of a linear growth promotion ($k \cdot P_s$) and a saturating Monod-type promotion ($k \cdot P_s/(K_p + P_s)$) are virtually indistinguishable. This creates a structural degeneracy where multiple models fit the noise floor equally well.

**2. The Autonomous Active DoE and LLM Intervention**

The Solution: To forcefully break this mathematical trap, the autonomous agent deliberately triggered an active Design of Experiments (DoE) pulse, injecting a substrate spike ($S_{sp} = 15.0$) at $t = 700$ hours. This forced the biological system into a highly dynamic transient state, rapidly driving up $P_s$ concentrations. To accurately evaluate these new transient dynamics, the LLM dynamically constructed and inserted the missing $dPs/dt$ ordinary differential equation directly into MATLAB, expanding the mathematical library in real time to capture the secondary metabolite's generation rate.

**3. The Winner: Saturating_Ps (Capturing Biological Limits)**

- *Physical Mechanism*: The Saturating Ps model utilizes a rational expression to describe how the secondary metabolite enhances or feeds back into the specific growth rate. Physically, this represents a metabolic pathway or cellular receptor system that can be stimulated by $P_s$, but only up to a finite maximum capacity ($k_{auto}$).

*Why it Won*: Once the DoE pulse pushed the system out of the low-concentration baseline, the true physical limits of the bioprocess were exposed. The Linear_Ps_Boost model wrongly assumed that biological stimulation scales infinitely, heavily overpredicting growth during the transient spike. Conversely, the Saturating_Ps model correctly captured the physiological plateau, mathematically restricting the growth promotion via its saturation constant ($K_p$). This structural accuracy yielded a near-perfect fit (MRV: 0.000297), The agent successfully mapped the correct mechanistic shape and locked in the final kinetic parameters ($k_{sat} = 0.1596$ and $K_{ps} = 6.55$), (true values are 0.15 and 5), proving the framework's ability to autonomously reconstruct complex physical truths from raw, hidden plant data.

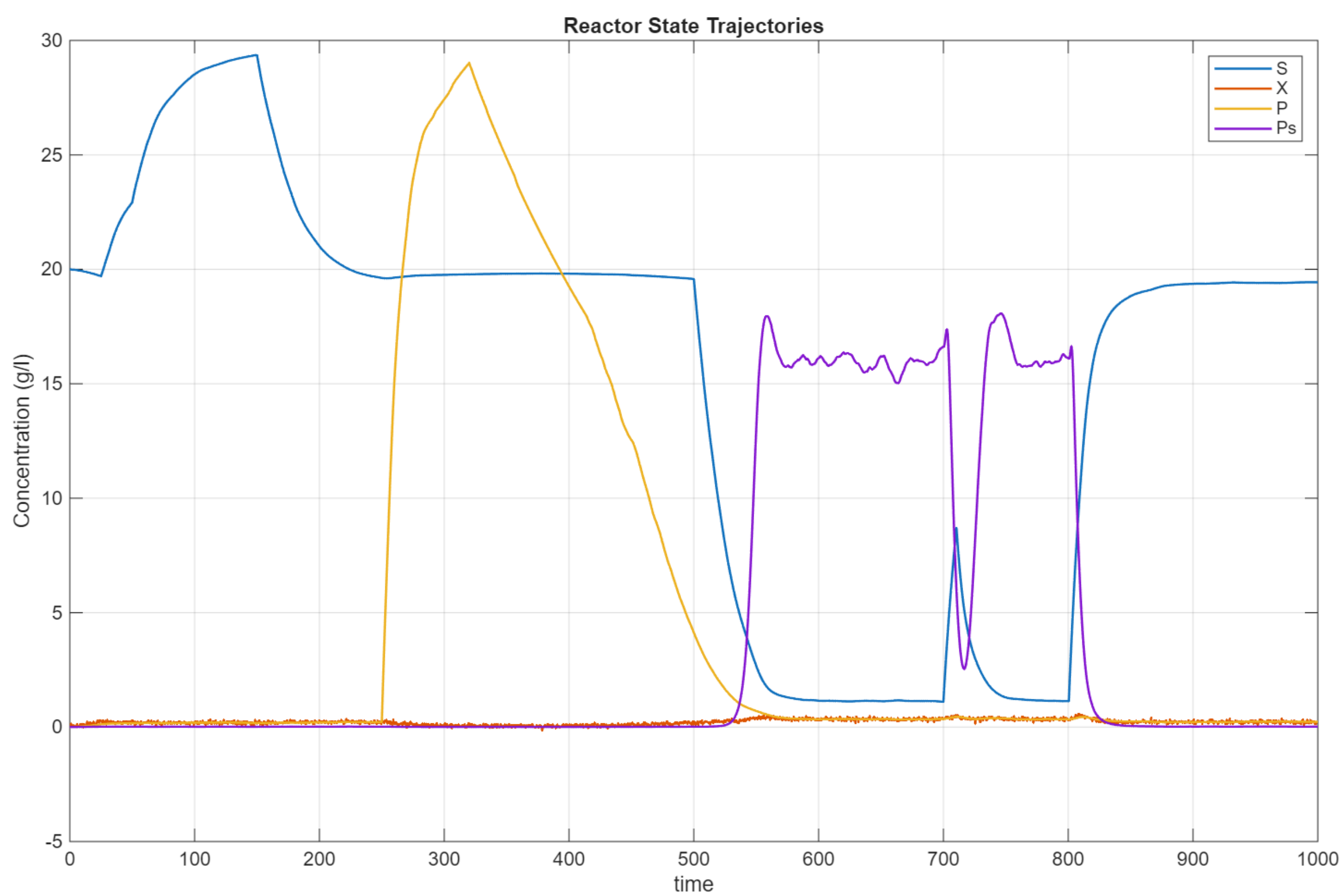

**Figure 6.** The operational regimes used in structural discovery.

The dynamic tracking performance of the MHE state estimator and the resulting isolation of the secondary additive error are illustrated in Figure 7. The temporal profiles compare the true underlying system growth rate ($\mu_{real}$), the trajectory captured by the Moving Horizon Estimation algorithm ($\mu_{MHE}$), and the structural prediction generated exclusively by the converged Level 1 multiplicative library ($\mu_{mult}$).

Throughout the primary operational window ($t < 500$ h), the $\mu_{mult}$trajectory aligns perfectly with the true system growth. This overlap mathematically validates the Stage 1 and Stage 2 tournaments, demonstrating that the combination of foundational Monod kinetics, Haldane substrate inhibition, and Luong product toxicity fully accounts for the macroscopic biological behavior under substrate-abundant conditions.

However, as the system transitions into the metabolic $P_s$ boost region ($t \geq 500$ h), a critical structural divergence emerges. While the MHE inherently relies on the system's differential mass balances to propagate the state, it operates dynamically by actively minimizing the objective function between the simulated states and the live sensor measurements over a rolling time window. As the primary substrate ($S$) depletes, the incomplete Level 1 structural prediction prematurely collapses toward zero (illustrated by the initial sharp decline of the blue dashed line, $\mu_{mult}$**)**. To mathematically reconcile the sustained biological activity observed by the sensors with this collapsing structural baseline, the MHE forces its estimated growth rate to capture the total physiological response, successfully maintaining $\mu_{MHE} \approx \mu_{real}$.

The macroscopic gap that forms between the MHE trajectory ($\mu_{MHE}$) and the baseline multiplicative prediction ($\mu_{mult}$) visually quantifies the unmodeled physics. Because the MHE successfully aggregates the total physiological response, the geometric difference between these two curves is mathematically equivalent to the missing additive flux ($\mu_{add}$). By dynamically subtracting the known $\mu_{mult}$baseline from the real-time $\mu_{MHE}$tracking data, the autonomous framework isolates this precise residual vector, allowing the agent to evaluate it against the Level 2 library and definitively resolve the saturating $P_s$mechanism during Stage 3.

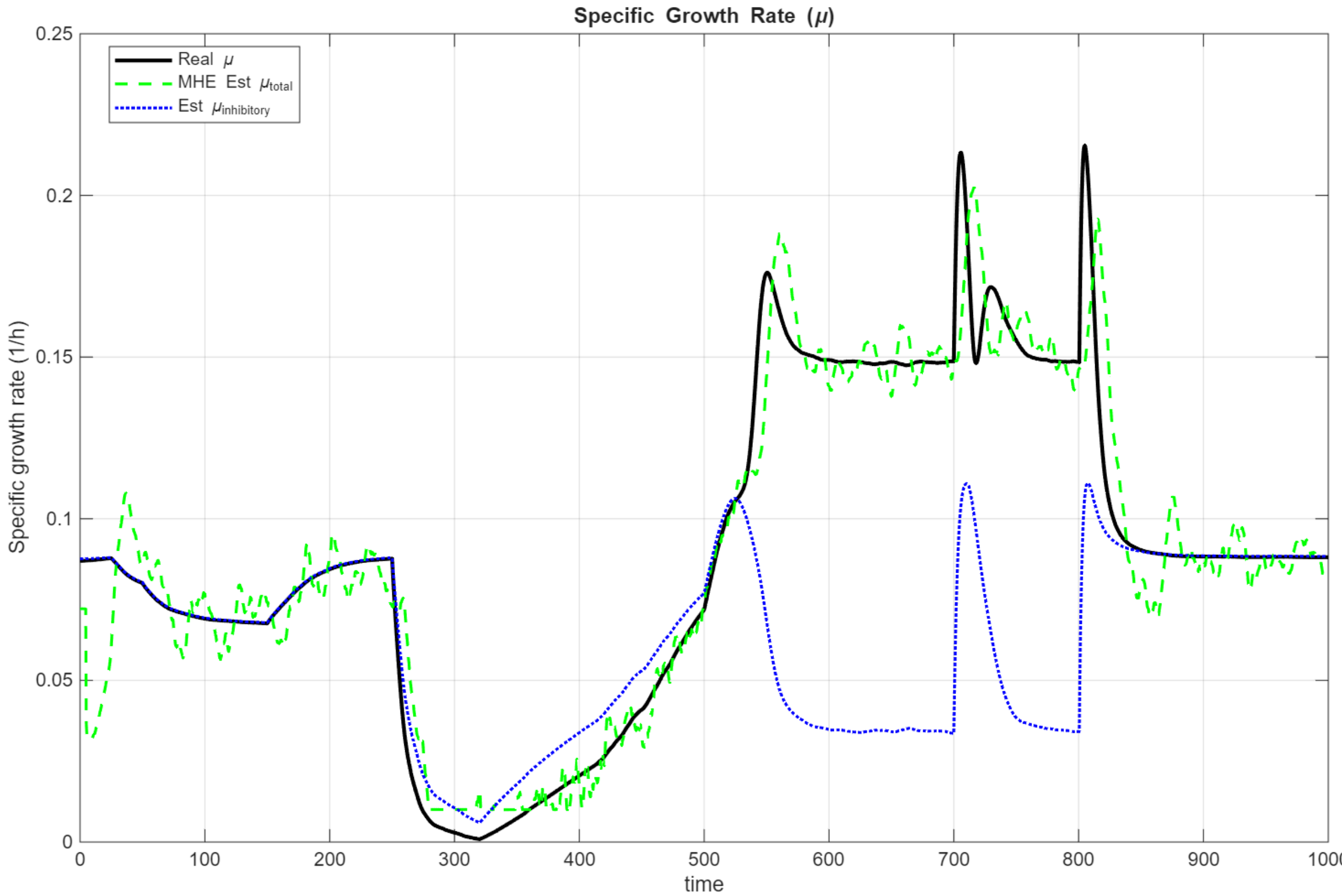

**Figure 7.** The specific growth rates.

Finally, the kinetic manifold alignments are shown in Figure 8. The discovery framework audits the empirical log residual ($E_{res}$) against the theoretical log and additive linear kinetic manifolds for each kinetic hypothesis. For the Haldane mechanism, the audit reveals a robust alignment between the scattered telemetry and the theoretical manifold within the inhibitory regime, confirming the structural validity of the kinetic basis. In contrast, the Luong toxicity audit demonstrates the framework's capacity to resolve non-linear manifold topologies, where the curvature of the innovation evidence strictly adheres to the power-law suppression signature defined by the toxicity manifold. The consistency between the scattered empirical evidence and the theoretical manifolds across both linear and non-linear regimes validates the agent's structural identification process, confirming that the discovered kinetics are the most parsimonious representations of the observed biological phenomena. The autocatalytic manifold audit reveals a saturating, concave-down topology, consistent with the identified Michaelis-Menten kinetic basis ($V_a = k_{auto}P_s/(K_p + P_s)$). Unlike the preceding phases which utilized logarithmic transformations to linearize multiplicative penalties, this audit is performed in the absolute-residual manifold, where the innovation evidence ($E_{res}$) is directly mapped to the additive growth rate contribution of the secondary metabolite. This curvature provides the primary evidence for the system's metabolic saturation, effectively rejecting linear boost distractors that would otherwise incorrectly predict unbounded growth. The alignment between the saturated evidence and the theoretical manifold—verified during the autonomous DoE intervention—validates that the secondary metabolite ($P_s$) acts as a rate-limited accelerator rather than a non-saturating feedback agent.

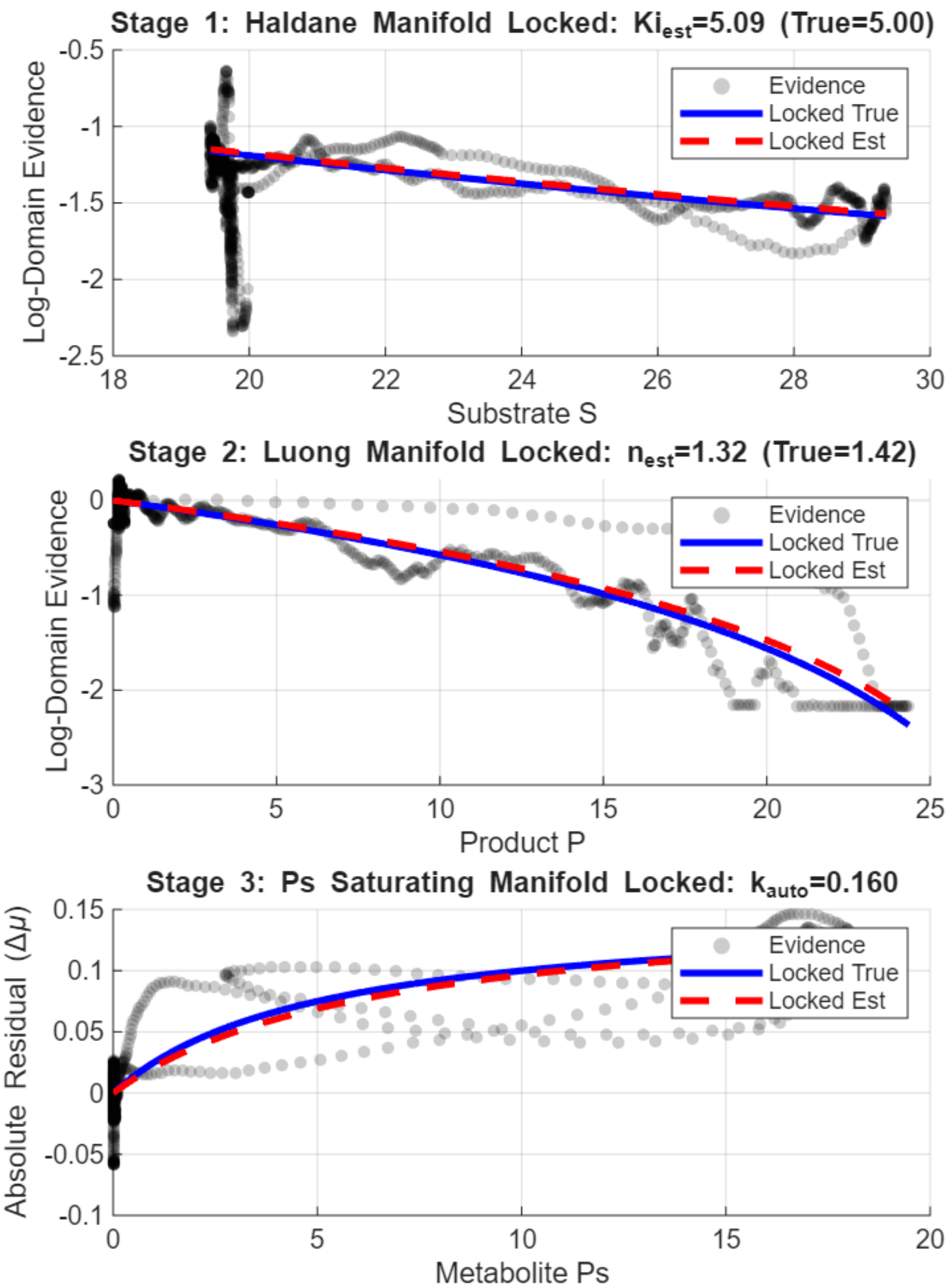

**Figure 8**. Discovery and reconstruction of multiplicative and additive kinetic manifolds.

### 4.4 Global Parameter Refinement

Following the sequential discovery of the kinetic manifold, the framework performed a global optimization to refine the complete parameter vector $\Theta_{global} = [\mu_{max}, K_s, K_i, P_m, n, K_{ib}, Y_{B/X}, k_{auto}, K_p]$. The objective function, $J_{global}(\Theta)$, was minimized across the entire 700-hour operational window using the previously discovered structural forms:

$$J_{global}(\Theta) = \sum_{k=1}^{N_{total}} \left(\mu_{obs}(k) - \hat{\mu}(\Theta, x(k))\right)^2 \quad (24)$$

The structural discovery architecture demonstrated high-fidelity parameter extraction, successfully resolving the system's kinetic topology through the sequential isolation of inhibitory and catalytic manifolds. As summarized in Table 2, the agent recovered the ground-truth kinetic parameters with high precision, achieving a recovery error of $< 6\%$ for the primary inhibitory and toxicity constants. The statistical audit confirms the robustness of the identified structure; the high $R^2$ values ($> 0.99$) and the negative AIC scores—indicating significant information gain—validate that the discovered Haldane-Luong-Autocatalytic model is the most parsimonious representation of the system's state-space telemetry. The narrow discrepancy between training and validation metrics ($SSE_{train} = 0.0434$ vs. $SSE_{val} = 0.0538$) further demonstrates that the framework avoids overfitting, maintaining structural integrity across distinct dynamic regimes.

**Table 2:** Global Parameter Identification and Statistical Audit

| Kinetic Parameter | Symbol | Identified Value | Ground Truth | Recovery Error (%) |
|---|---|---|---|---|
| Haldane Inhibition | $K_i$ | **5.02188** | 5.0000 | +0.44% |
| Toxicity Limit | $P_m$ | **29.1167** | 30.0000 | -2.94% |
| Toxicity Exponent | $n$ | **1.3351** | 1.4200 | -5.98% |
| Autocatalytic Rate | $k_{auto}$ | **0.1472** | 0.1500 | -1.87% |

| **Kinetic Parameter** | **Symbol** | **Identified Value** | **Ground Truth** | **Recovery Error (%)** |
| --- | --- | --- | --- | --- |
| Saturation Constant | $K_p$ | **5.5065** | 5.0000 | +10.13% |

**Statistical Model Verification**

| **Metric** | **Training Set** | **Validation Set** |
| --- | --- | --- |
| Sum of Squared Errors (SSE) | 0.0434 | 0.0538 |
| Akaike Information Criterion (AIC) | $-14{,}535.83$ | $-14{,}233.34$ |
| Coefficient of Determination ($R^2$) | 0.9965 | 0.9917 |

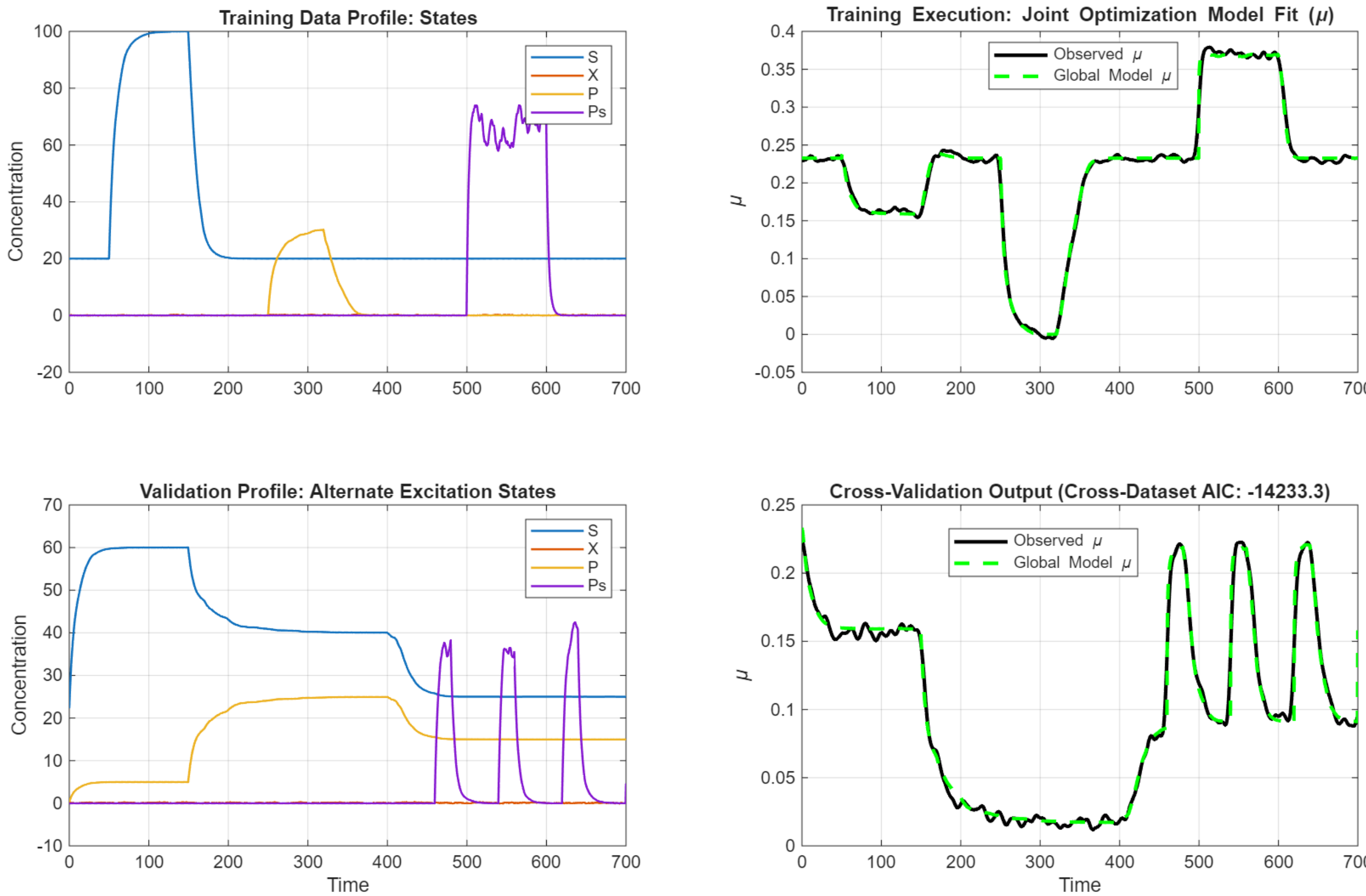

**Figure 9.** Final validation and global parameter estimation.

## 5. LIMITATIONS AND FUTURE WORK

While the hybrid autonomous discovery framework demonstrates a high capacity for resolving complex structural kinetic ambiguities, transitioning from an offline-assisted design to a fully online, closed-loop LLM controller presents several significant technical challenges. This study deliberately adopted a hybrid approach—where the AI conducts symbolic reasoning and experimental design offline to construct scheduled interventions—due to three primary limitations of integrating current LLM infrastructure with high-frequency numerical computing environments.

### 5.1. Computational Latency and API Bottlenecks

The real-time execution of bioprocess digital twins requires solving differential equations at high frequencies. Integrating live HTTP requests to external LLM application programming interfaces (APIs) introduces severe asynchronous latency. A single inference call to evaluate $\Delta$AIC and propose an actuator command can take several seconds to process. When placed inside an active process control and state estimation loop, this latency can create substantial computational bottleneck, ultimately destabilizing controller timing and real-time state tracking.

### 5.2. Non-Deterministic Output and Execution Fragility

Industrial process control demands strict deterministic reliability. Standard large language models are inherently probabilistic and prone to structural formatting deviations. Even with highly constrained prompting, online APIs can occasionally generate responses with conversational text outside of the requested JSON schema, or output unparsable data types. Deploying a live LLM directly to physical actuators without extensive digital safety wrappers risks unhandled code exceptions or, in a physical deployment, catastrophic reactor failure resulting from out-of-bounds setpoint commands.

### 5.3. Scientific Reproducibility

A core requirement of computational kinetic research is reproducibility. If an active Design of Experiments (DoE) pulse is generated by a live, online API call during simulation, the inherent variance in LLM responses dictates that subsequent executions may yield slightly different setpoint commands or trigger times. By executing the agent's symbolic reasoning offline and formally scheduling the resulting interventions into the operational baseline, the methodology presented here guarantees mathematical consistency and allows independent researchers to perfectly reproduce the structural resolution of the identifiability traps.

**Future Directions**

Future iterations of this framework will focus on bridging the gap between theoretical autonomy and practical real-time execution. A primary avenue of research involves transitioning away from cloud-based LLM APIs in favor of deploying quantized, domain-specific small language models (SLMs) locally on the edge. Embedding lightweight, specialized models directly within the numerical execution environment would eliminate network latency and drastically reduce parsing errors. Furthermore, developing deterministic software guardrails—symbolic filtering layers that intercept, mathematically verify, and strictly bound LLM actuator commands prior to execution—will be critical for safely advancing this architecture from simulated digital twins to fully closed-loop active experimentation in physical wet-lab environments. A concrete instantiation of this principle has recently been proposed by Nogueira and Skogestad (2026), who map each loop of a classical Advanced Regulatory Control (ARC) chain onto a specialised LLM operator agent and encapsulate all MIN/MAX selector and override logic within a single deterministic orchestrator; because every constraint conflict is resolved by this orchestrator rather than by the LLM itself, the resulting multi-agent system inherits the safety guarantees of the underlying control chain regardless of the language model's output. Adopting an analogous orchestrator-mediated architecture, in which the LLM Discovery Agent proposes candidate structural hypotheses but a deterministic supervisory layer enforces physical bounds and arbitrates conflicting model recommendations, offers a promising path toward safely closing the loop between structural discovery and online actuation in this framework.

A fundamental limitation of the current framework is its reliance on a Closed-World RAG architecture. The discovery engine is bounded by the Level 1 and Level 2 candidate libraries curated offline by the human supervisor. As demonstrated by the additive illusion, if the true governing physics are absent from the mathematical vocabulary, the agent is forced into constructing a mathematical fiction out of the available terms. To achieve true, unconstrained autonomy, future iterations of this architecture must transition from a static library dependency to open-ended generative discovery.

**5.4. Open-World RAG and Dynamic Library Expansion** Rather than relying on offline curation, the LLM agent can be equipped with external tool-calling capabilities (e.g., Python-based web crawlers and API integrations with academic databases). In an Open-World RAG paradigm, if the agent detects a persistent Identifiability Trap or fails to lower the Mean Residual Variance below the noise floor using the existing JSON repository, it would autonomously query literature databases. By scraping recent biochemical publications for novel kinetic models, parsing the text into mathematical structures, and dynamically injecting these new equations into the JSON blueprint, the agent could self-expand its mathematical vocabulary in real-time.

**5.5. Symbolic Regression and Neural Fallbacks** In scenarios where the missing physics represents a truly novel phenomenon not yet documented in the literature, RAG retrieval will inherently fail. To prevent the framework from forcing an incorrect parametric fit, the agent's logic must include an ultimate fallback mechanism. Upon exhausting the expanded RAG libraries, the agent would transition from model selection to equation generation. This can be achieved through LLM-guided Symbolic Regression—where the agent mathematically constructs novel structural equations from base algebraic operators—or by bridging into Scientific Machine Learning (SciML). In the latter case, the agent would instruct the numerical solver to embed a Multi-Layer Perceptron (MLP) directly into the differential equations, forming a Universal Differential Equation (UDE). The neural network would universally approximate the unknown residual flux, mapping the blind spots of the mechanistic library without falsely constraining them to known kinetic forms. This UDE-based fallback is not merely speculative: Nogueira and coworkers have already demonstrated systematic hybrid-model and governing-equation discovery for non-linear transport and sorption systems using precisely this scientific machine learning strategy (Santana et al., 2023), and their approach offers a natural template for embedding a neural residual term directly into the CSTR mass balances whenever the RAG libraries are exhausted.

## 6. CONCLUSION

This work demonstrates an autonomous, LLM-driven diagnostic tool capable of overcoming structural model mismatches in complex bioprocesses. By combining MHE state estimation, logarithmic manifold decomposition, and a Two-Level RAG library, the framework systematically isolates and identifies hidden multiplicative and additive kinetic terms from dynamic telemetry. Crucially, the integration of an Active DoE controller enables the agent to break mathematical identifiability traps through online experimental forcing. While currently bounded by closed-world RAG libraries and computational latency, future extensions toward local edge language models, open-world literature retrieval, and symbolic/neural fallbacks offer a clear path toward fully autonomous physical lab deployments.

## 7. REFERENCES

1. Aiba, S., Shoda, M., & Nagatani, M. (1968). Kinetic study of ethanol inhibition in fermentation. *Biotechnology and Bioengineering*, 10(6), 845–864.
2. Andrews, J. F. (1968). A mathematical model for the continuous culture of microorganisms utilizing inhibitory substrates. *Biotechnology and Bioengineering*, 10(6), 707–723.
3. Boiko, D. A., MacKnight, R., Kline, B., & Gomes, G. (2023). Autonomous chemical research with large language models. *Nature*, 624(7992), 570–578.
4. Bran, A. M., Cox, S., Schilter, O., Baldassari, R., White, A. D., & Jablonka, K. M. (2024). Augmenting large language models with chemistry tools. *Nature Machine Intelligence*, 6(5), 525–535.
5. Brunton, S. L., Proctor, J. L., & Kutz, J. N. (2016). Discovering governing equations from data by sparse identification of nonlinear dynamical systems. *Proceedings of the National Academy of Sciences*, 113(15), 3932–3937.
6. Chang, A., Scheer, M., Grote, A., Schomburg, I., & Schomburg, D. (2008). BRENDA, AMNDA and FRENDA the enzyme information system: new content and tools in 2008. *Nucleic Acids Research*, 36(suppl_1), D280–D284.
7. de Carvalho Servia, A., Mowbray, M., & del Rio-Chanona, E. A. (2024). Automated dynamic kinetic model discovery via information-theoretic genetic programming.

8. Herbert, D. (1958). Some principles of continuous culture. *Recent Progress in Microbiology*, 381–396.
9. Hill, A. V. (1910). The possible effects of the aggregation of the molecules of hæmoglobin on its dissociation curves. *The Journal of Physiology*, 40(4), iv–vii.
10. Hinshelwood, C. N. (1946). *The Chemical Kinetics of the Bacterial Cell*. Oxford University Press.
11. Jerusalimsky, N. D. (1967). Bottlenecks in physiological processes. *Microbial Physiology and Continuous Culture*, 23–32.
12. LeCun, Y., Bengio, Y., & Hinton, G. (2015). Deep learning. *Nature*, 521(7553), 436–444.
13. Levenspiel, O. (1980). *The Chemical Reactor Omnibook*. OSU Book Stores.
14. Lewis, P., Perez, E., Pichtulis, A., Petroni, F., Karpukhin, V., Goyal, N., ... & Riedel, S. (2020). Retrieval-augmented generation for knowledge-intensive NLP tasks. *Advances in Neural Information Processing Systems*, 33, 9459–9474.
15. Luong, J. H. T. (1985). Kinetics of ethanol inhibition in alcohol fermentation. *Biotechnology and Bioengineering*, 27(3), 280–285.
16. Luong, J. H. T. (1987). Generalization of Monod kinetics for toxicity of substrate and product in bioreactors and application to ethanol fermentation. *Biotechnology and Bioengineering*, 31(3), 242–247.
17. Monod, J. (1949). The growth of bacterial cultures. *Annual Review of Microbiology*, 3(1), 371–394.
18. Mowbray, M., Savage, T., Wu, C., Zhang, D., & del Rio-Chanona, E. A. (2022). Automated kinetic model discovery for bioprocesses. *Computers & Chemical Engineering*.
19. Nogueira, R. A., & Skogestad, S. (2026). Orchestrator-mediated multi-agent LLM control chains for regulatory control.
20. Pajak, M., et al. (2025). Multi-agent LLM frameworks for operational decision-making, flowsheet synthesis, and simulation workflows in chemical process contexts.
21. Raissi, M., Perdikaris, P., & Karniadakis, G. E. (2019). Physics-informed neural networks: A deep learning framework for solving forward and inverse problems involving nonlinear partial differential equations. *Journal of Computational Physics*, 378, 686–707.
22. Rawlings, J. B., & Mayne, D. Q. (2009). *Model Predictive Control: Theory and Design*. Nob Hill Publishing.

23. Santana, et al. (2023). Governing equation discovery for nonlinear transport and sorption systems using scientific machine learning.
24. Schmidt, A., et al. (2022). Continuous biomanufacturing and process analytical technology initiatives in the biopharmaceutical industry.
25. Webb, F. C. (1963). *Biochemical Engineering*. D. Van Nostrand Company.
26. Yano, T., & Koga, S. (1969). Dynamic behavior of the chemostat subject to substrate inhibition. *Biotechnology and Bioengineering*, 11(2), 139–153.